\documentclass[11pt]{amsart} 
\usepackage{amscd,amssymb,amsxtra}
\usepackage{mathrsfs}  
\DeclareSymbolFontAlphabet{\mathrsfs}{rsfs}
\usepackage[mathscr]{eucal} 
\usepackage{mathabx}\DeclareMathSymbol{\emptyset}{\mathord}{symbols}{59}
                    \DeclareMathSymbol{\boxtimes}{\mathbin}{AMSa}{"02}
\usepackage{comment}
\usepackage{color}
\usepackage{enumitem} 
\usepackage{marginnote}

\makeatletter
\let\@secnumfont\bfseries
\def\section{\@startsection{section}{1}%
  \z@{4\linespacing\@plus\linespacing}{\linespacing}%
  {\bfseries\centering}}
\def\introsection{\@startsection{section}{1}%
  \z@{3\linespacing\@plus\linespacing}{\linespacing}%
  {\bfseries\centering}}
\def\subsection{\@startsection{subsection}{2}%
   \z@{1.25\linespacing\@plus.7\linespacing}{.5\linespacing}%
   {\normalfont\bfseries}}
\def\subsectionsinline{\def\subsection{\@startsection{subsection}{2}%
  \z@{1\linespacing\@plus.7\linespacing}{-.5em}%
  {\normalfont\bfseries}}}

\makeatother

\theoremstyle{definition}
\newtheorem{definition}[equation]{Definition}
\newtheorem{example}[equation]{Example}

\newtheorem*{definition*}{Definition}
\newtheorem*{example*}{Example}
\newtheorem*{problem*}{\color{blue}Problem}
\newtheorem*{exercise*}{Exercise}
\newtheorem*{question*}{\color{blue}Question}
\newtheorem*{project*}{\color{blue}Project}
\newtheorem*{construction*}{Construction}

\theoremstyle{remark}

\newtheorem{remark}[equation]{Remark}

\newtheorem*{note*}{Note}
\newtheorem*{notation*}{Notation}
\newtheorem*{remark*}{Remark}
\newtheorem*{data*}{Data}

\theoremstyle{plain}
\newtheorem{theorem}[equation]{Theorem}
\newtheorem{corollary}[equation]{Corollary}

\newtheorem*{theorem*}{Theorem}
\newtheorem*{corollary*}{Corollary}
\newtheorem*{lemma*}{Lemma}
\newtheorem*{proposition*}{Proposition}
\newtheorem*{conjecture*}{Conjecture}
\newtheorem*{claim*}{Claim}
\newtheorem*{proposal*}{Proposal}
\newtheorem*{conclusion*}{Conclusion}
\newtheorem*{hypothesis*}{Hypothesis}
\newtheorem*{assumption*}{Assumption}

\newenvironment{proof*}[1][\proofname]{
  \begin{proof}[#1]}{  
\end{proof}}

\numberwithin{equation}{section}

\definecolor{refkey}{rgb}{0,.6,.4}

\renewcommand{\:}{\colon}
\renewcommand{\AA}{{\mathbb A}}
\newcommand{\Ahat}{{\hat A}}
\DeclareMathOperator{\Aut}{Aut}
\newcommand{\CC}{{\mathbb C}}
\newcommand{\CP}{{\mathbb C\mathbb P}}

\newcommand{\EE}{\mathbb E}

\DeclareMathOperator{\End}{End}

\newcommand{\HH}{{\mathbb H}}
\DeclareMathOperator{\Hom}{Hom}

\DeclareMathOperator{\Met}{Met}

\newcommand{\PP}{{\mathbb P}}
\DeclareMathOperator{\pt}{pt}
\newcommand{\QQ}{{\mathbb Q}}

\newcommand{\RR}{{\mathbb R}}
\newcommand{\TT}{\mathbb T}
\DeclareMathOperator{\Spin}{Spin}

\DeclareMathOperator{\Tr}{Tr}
\DeclareMathOperator{\tr}{tr}

\newcommand{\ZZ}{{\mathbb Z}}
\DeclareMathOperator{\ch}{ch}
\newcommand{\chiup}{\raise.5ex\hbox{$\chi$}}
\newcommand{\cir}{S^1}
\DeclareMathOperator{\coker}{coker}
\newcommand{\dbar}{{\overline\partial}}

\DeclareMathOperator{\ind}{ind}
\newcommand{\inv}{^{-1}}
\DeclareRobustCommand{\mstrut}{^{\vphantom{1*\prime y\vee M}}}

\DeclareMathOperator{\rank}{rank}
\newcommand{\res}[1]{\negmedspace\bigm|\mstrut_{#1}}

\newcommand{\temsquare}{\raise3.5pt\hbox{\boxed{ }}}

\newcommand{\zmod}[1]{\ZZ/#1\ZZ}

\newcommand{\zt}{\zmod2}

\DeclareMathOperator{\SO}{SO}

\let\O\relax
\DeclareMathOperator{\O}{O}

\DeclareMathOperator{\PSL}{PSL}

\usepackage[all,2cell]{xy}\renewcommand{\cir}{\ensuremath{S^1}}
\usepackage[colorlinks,backref=page,citecolor=refkey]{hyperref}
\let\O\relax\DeclareMathOperator{\O}{O}
\DeclareMathOperator{\Bord}{Bord}
\DeclareMathOperator{\Cliff}{Cliff}
\DeclareMathOperator{\Det}{Det}
\DeclareMathOperator{\Ell}{Ell}
\DeclareMathOperator{\Euler}{Euler}
\DeclareMathOperator{\Fix}{Fix}
\DeclareMathOperator{\Fred}{Fred}
\DeclareMathOperator{\Imm}{Im}
\DeclareMathOperator{\Line}{Line}
\DeclareMathOperator{\Sign}{Sign}
\DeclareMathOperator{\Sym}{Sym}
\DeclareMathOperator{\Todd}{Todd}
\DeclareMathOperator{\Vect}{Vect}
\DeclareMathOperator{\Vol}{Vol}
\DeclareMathOperator{\ad}{ad}
\DeclareMathOperator{\curv}{curv}
\DeclareMathOperator{\hol}{hol}
\DeclareMathOperator{\image}{image}
\DeclareMathOperator{\ord}{ord}
\DeclareMathOperator{\pfaff}{pfaff}
\DeclareMathOperator{\pr}{pr}
\DeclareMathOperator{\sign}{sign}
\DeclareMathOperator{\spec}{spec}
\newcommand{\AP}{\ssA_{\sP}}
\newcommand{\BSO}{B\!\SO}
\newcommand{\Bd}{\Bord_{\langle d,d+1  \rangle}}
\newcommand{\Clpm}{\Cliff_{\pm n}}
\newcommand{\Cnm}{\Cliff_{-n}}
\newcommand{\Cnp}{\Cliff_{+n}}
\newcommand{\Cn}{\Cliff_{n}}
\newcommand{\Fn}{\Fred_n}
\newcommand{\Mn}{\mathbf{Man}_n}
\newcommand{\PH}{\PP\sH}
\newcommand{\RZ}{\RR/\ZZ}
\newcommand{\Rn}{\RR^n}
\newcommand{\Set}{\mathbf{Set}}
\newcommand{\Vk}{\Vect_k}
\newcommand{\bX}{\partial X}
\newcommand{\cP}{\sC_{\sP}}
\newcommand{\et}[1]{\eta \mstrut _{#1}}
\newcommand{\evo}[2]{e^{-i(t_{#1}-t_{#2})H/\hbar}}
\newcommand{\hooklongrightarrow}{\lhook\joinrel\longrightarrow} 
\newcommand{\op}{^{\textnormal op}}
\newcommand{\q}{^{(q)}}
\newcommand{\sA}{\mathcal{A}}
\newcommand{\sB}{\mathscr{B}}
\newcommand{\sC}{\mathscr{C}}
\newcommand{\sE}{\mathscr{E}}
\newcommand{\sG}{\mathscr{G}}
\newcommand{\sH}{\mathscr{H}}
\newcommand{\sL}{\mathcal{L}}
\newcommand{\sM}{\mathscr{M}}
\newcommand{\sO}{\mathcal{O}}
\newcommand{\sP}{\mathscr{P}}
\newcommand{\sS}{\mathscr{S}}
\newcommand{\scrF}{\mathscr{F}}
\newcommand{\sect}[1]{\sL(#1)}
\newcommand{\ssA}{\mathfrak{A}}
\newcommand{\tD}{\widetilde{D}}
\newcommand{\tE}{\widetilde{E}}

\newcommand{\tX}{\widetilde{X}}

\newcommand{\tsH}{\widetilde{\sH}}
\newcommand{\xia}[1]{\xi \mstrut _{#1}}
\newcommand{\zD}{\zeta \mstrut _\Delta }
\renewcommand{\Im}{\Imm}
\renewcommand{\Re}{\text{Re}}
\renewcommand{\SS}{\mathbb{S}}

\begin{document}

\abovedisplayskip18pt plus4.5pt minus9pt
\belowdisplayskip \abovedisplayskip
\abovedisplayshortskip0pt plus4.5pt
\belowdisplayshortskip10.5pt plus4.5pt minus6pt
\baselineskip=15 truept
\marginparwidth=55pt

\makeatletter
\renewcommand{\tocsection}[3]{%
  \indentlabel{\@ifempty{#2}{\hskip1.5em}{\ignorespaces#1 #2.\;\;}}#3}
\renewcommand{\tocsubsection}[3]{%
  \indentlabel{\@ifempty{#2}{\hskip 2.5em}{\hskip 2.5em\ignorespaces#1%
    #2.\;\;}}#3} 
\renewcommand{\tocsubsubsection}[3]{%
  \indentlabel{\@ifempty{#2}{\hskip 5.5em}{\hskip 5.5em\ignorespaces#1%
    #2.\;\;}}#3} 
\def\@makefnmark{%
  \leavevmode
  \raise.9ex\hbox{\fontsize\sf@size\z@\normalfont\tiny\@thefnmark}} 
\def\multfoot{\textsuperscript{\tiny\color{red},}}
\def\footref#1{$\textsuperscript{\tiny\ref{#1}}$}
\makeatother

\setcounter{tocdepth}{2}


 \title[Atiyah-Singer Index Theorem]{The Atiyah-Singer Index Theorem} 
 \author[D. S. Freed]{Daniel S.~Freed}
 \thanks{This material is based upon work supported by the National Science
Foundation under Grant Number DMS-2005286.}
 \address{Department of Mathematics \\ University of Texas \\ Austin, TX
78712} 
 \email{dafr@math.utexas.edu}
 \dedicatory{In memory of Michael Atiyah}
 \date{July 7, 2021}
 \begin{abstract} 
 The Atiyah-Singer index theorem, a landmark achievement of the early 1960s,
brings together ideas in analysis, geometry, and topology.  We recount some
antecedents and motivations; various forms of the theorem; and some of its
implications, which extend to the present.
 \end{abstract}
\maketitle

{\small
\def\reftext{References}
\renewcommand{\tocsection}[3]{%
  \begingroup 
   \def\tmp{#3}%
   \ifx\tmp\reftext
  \indentlabel{\phantom{1}\;\;} #3%
  \else\indentlabel{\ignorespaces#1 #2.\;\;}#3%
  \fi\endgroup}
\tableofcontents
}

   \section{Introduction}\label{sec:1}

Consider the Riemann sphere $\CP^1=\CC\cup\{\infty \}$.  Let
$\{z_i\}_{i=1}^N\subset \CP^1$ be a finite set, and to each $i\in \{1,\dots
,N\}$ suppose a nonzero integer~$n_i$ is given.  A classical problem asks for
a meromorphic function~$f$ with a zero or pole at each~$z_i$.  If $n_i>0$,
then $z_i$~is a zero of multiplicity~$n_i$; if $n_i<0$, then $z_i$~is a pole
of order~$|n_i|$.  The solution is straightforward.  Namely,
$f(z)=c\prod(z-z_i)^{n_i}$ is a rational function, unique up
to the constant~$c\in \CC\setminus \{0\}$.  In other words, now
allowing~$f\equiv 0$, the solutions form a one-dimensional complex vector
space.  If we replace~$\CP^1$ by a closed Riemann surface of positive genus,
then there is an obstruction to existence of a meromorphic function with
specified zeros and poles.  For example, an elliptic curve can be realized as
a quotient $E=\CC/(\ZZ+\ZZ\tau )$ of the complex line by the full lattice
generated by~$1,\tau $ for some~$\tau \in \CC$ with $\Im\tau >0$.  A
meromorphic function on~$E$ lifts to a doubly periodic function on~$\CC$, and
the single constraint on the zeros and poles of a meromorphic function is
$\sum n_iz_i\in \ZZ+\ZZ\tau $.  Proceeding from~$E$ to a
general closed Riemann surface~$X$, we encounter more constraints.  In fact,
the constraints form a vector space whose dimension is the genus of~$X$, a
topological invariant.  Meromorphic functions are solutions to the
Cauchy-Riemann equation, a linear elliptic partial differential equation.
The solutions and obstructions to this elliptic PDE are ``counted'' via
topology.  In a more general form, this is the classical Riemann-Roch
theorem.
 
The Atiyah-Singer index theorem, formulated and proved in 1962--3, is a vast
generalization to arbitrary elliptic operators on compact manifolds of
arbitrary dimension.  The Fredholm index in question is the dimension of the
kernel minus the dimension of the cokernel of a linear elliptic operator.
The Atiyah-Singer theorem computes the index in terms of topological
invariants of the operator and topological invariants of the underlying
manifold.  The theorem weaves together concepts and results in algebraic
topology, algebraic geometry, differential geometry, and linear analysis; its
ramifications go far beyond, in number theory, representation theory,
operator algebras, nonlinear analysis, and theoretical physics.  Furthermore,
index theory itself is a sprawling enterprise.  The basic Atiyah-Singer
theorem spawned numerous generalizations and novel pathways.  This paper---a
tribute to Michael Atiyah---naturally focuses on aspects of his work and his
influence.  Even thus restricted, we can only skim the surface of this rich
story.
 
There are antecedents of the index theorem from algebraic geometry and
topology on the one hand, and from analysis on the other.  We discuss these
in turn in~\S\ref{sec:2} and~\S\ref{sec:3}.  The basic Atiyah-Singer theorem
is the subject of~\S\ref{sec:4}.  The first proof is based on cobordism and
in broad outline follows Hirzebruch's proofs of his signature and
Riemann-Roch theorems.  The second proof is based on $K$-theory; it is
inspired by Grothendieck's Riemann-Roch theorem.  In~\S\ref{sec:5} we take up
some of the extensions and variations of the basic theorem.  These include an
equivariant index theorem, the index theorem for parametrized families of
operators, the index theorem for manifolds with boundary, and a few more.  At
this point our exposition makes a transition from global topological
invariants of general linear elliptic operators to local geometric invariants
of Dirac operators.  Heat equation methods are the subject of~\S\ref{sec:6},
the first application being a local index theorem.  New geometric invariants
of Dirac operators appear in~\S\ref{sec:7}.  In~\S\ref{sec:8} we turn to
physics, which was Atiyah's focus after the mid 1980's and which provided an
unanticipated playground for the circle of ideas surrounding the index
theorem.  We focus on anomalies in quantum theory, a subject to which Atiyah
and Singer made an early contribution.  

Each section of this paper has more introductory material, which we recommend
even to the casual reader.  Also, a lecture based on this paper may be viewed
at~\cite{F1}.

Michael had great mathematical and personal charisma.  His writings capture
his vibrancy, as did his lectures, some of which are available online.  He
wrote many wonderful expository articles about the index theorem, especially
of the early period; you will enjoy perusing them.  

I warmly thank Simon Donaldson, Charlie Reid, and Graeme Segal for their
careful reading of and comments on an earlier version.

   \section{Antecedents and motivations from algebraic geometry and topology}\label{sec:2}

Enumerative problems in algebraic geometry often lead to integers that have a
topological interpretation.  A classical example is the Riemann-Roch formula,
which is our starting point in~\S\ref{subsec:2.1}.  The higher dimensional
generalization was taken up by Fritz Hirzebruch in the early 1950's, as we
recount in~\S\ref{subsec:2.3}.  A few years later Alexander Grothendieck
extended Hirzebruch's theorem to a relative version, that is, to proper maps
of complex manifolds.  In the process he introduced $K$-theory for sheaves.
His ideas, briefly presented in~\S\ref{subsec:2.2}, play a fundamental role
in variations of the Atiyah-Singer index theorem a decade later.
More immediately, as Graeme Segal writes in this volume~\cite{Seg4}, Atiyah
and Hirzebruch transported Grothendieck's $K$-theory over to algebraic
topology.  Raoul Bott's computation of the stable homotopy groups of Lie
groups, which took place during the same period as Hirzebruch's and
Grothendieck's work on the Riemann-Roch theorem, is the cornerstone of their
theory.  Crucial for the index theorem are the resulting integrality
theorems, of which we mention a few in~ \S\ref{subsec:2.4}.  This led to a
question---Why is the $\Ahat$-genus an integer for a spin manifold?---which
in early 1962 was the immediate catalyst for Atiyah and Singer's
collaboration.

  \subsection{The Riemann-Roch theorem}\label{subsec:2.1}

 Let $X$ be a smooth projective curve over~$\CC$, i.e., a one-dimensional
closed complex submanifold of a complex projective space.  A divisor~$D$ is a
finite set of points on~$X$ with an integer~$\ord_x(D)$ attached to each
point~$x\in D$.  A divisor determines a holomorphic line bundle on~$X$; let
$\sect D$~denote the space of holomorphic sections of this bundle.  We can
describe~$\sect D$ as the space of meromorphic functions on~$X$ which have a
pole of order~ $\le \ord_x(D)$ at each~$x\in X$.  A basic problem in the
theory of curves is: Compute the dimension of~$\sect D$.  While this is quite
difficult in general, there is a {\it topological\/} formula for $\dim\sect D
- \dim\sect{K-D}$, where $K$~is a \emph{canonical divisor} of~$X$.  (The zero
set of a holomorphic 1-form, weighted by the orders of the zeros, is a
canonical divisor.)  

  \begin{theorem}[Riemann-Roch]\label{thm:2}
 Let $X$~be a smooth projective curve and let $D$~be a divisor on~$X$.  Then
  \begin{equation}\label{eq:1}
     \dim\sect D - \dim\sect{K-D} = \deg(D) -g + 1. 
  \end{equation}
  \end{theorem}

\noindent 
 Here $g$~is the genus of the curve~$X$, its fundamental topological
invariant, which is defined to be $\frac{1}{2}\rank H^1(X;\ZZ)$.  Also,
$\deg(D)=\sum \ord_x(D)$~is the sum of the integers which define the
divisor~$D$.  If $\deg(D) > 2g-2$, it can be shown that $\sect{K-D}=0$, so
that in that case \eqref{eq:1}~ provides a complete solution to the problem
of computing $\dim\sect D$.  Theorem~\ref{thm:2} is the classical
Riemann-Roch\footnote{Riemann~\cite{Ri} proved the inequality $\dim\sect D
\ge \deg(D) -g + 1$, and then Roch~\cite{Ro} proved the more
precise~\eqref{eq:1}.  Sadly, Roch died of tuberculosis at the age of~26,
just months after the 39~year-old Riemann succumbed to tuberculosis.}
formula.  The Atiyah-Singer index theorem is a vast generalization
of~\eqref{eq:1}, as we will see.

Let us immediately note one consequence of the Riemann-Roch formula.  Take
$D=\sO$ to be the trivial divisor consisting of no points.  Then $\sect
{\sO}$~is the space of constant functions and $\sect K$~is the space of
holomorphic differentials.  We deduce from~\eqref{eq:1} that the latter has
dimension~$g$.  It follows that $g$~is an integer, i.e., $\rank H^1(X;\ZZ)$
is even.  Therefore, one-half the Euler number~$\Euler(X)$ is an integer, our
first example of an {\it integrality theorem}.  The proof is noteworthy:
$1-\Euler(X)/2$~is an integer because it is the dimension of a vector space,
namely~$\sL(K)$.
 
In the last decade of the $19^{\textnormal{th}}$~century, Noether, Enriques,
and Castelnuovo generalized the Riemann-Roch inequality and equality to
algebraic surfaces; see~\eqref{eq:4} below.

  \subsection{Hirzebruch's Riemann-Roch and Signature
Theorems}\label{subsec:2.3} 

We skip far ahead to the years 1945--1954 and the work of young Hirzebruch,
based on two important developments in geometry.  The first, initiated by
Leray, is the theory of sheaves.  The second are the results in Thom's
thesis, particularly those concerning bordism\footnote{Thom, Hirzebruch, and
many others use `cobordism' in place of `bordism'; Atiyah~\cite{A10}
clarified the relationship.} groups of smooth manifolds.
We state two of Hirzebruch's main results, which are recounted in~\cite{H1}.

Let $X$~be a nonsingular projective variety of complex dimension~$n$ and
$V\to X$~a holomorphic vector bundle.  (In our discussion of curves we used
divisors; recall that a divisor determines a holomorphic line bundle, which
makes the link to our formulation here.)  Then the cohomology
groups~$H^q(X,V)$ are defined via sheaf theory: $H^0(X,V)$~is the vector
space of holomorphic sections of $V\to X$, and $H^q(X,V)$ for $q\ge1$ are
derived from resolutions of the sheaf of holomorphic sections of $V\to X$.
The cohomology groups are finite dimensional, which can be proved using the
theory of elliptic differential operators and Dolbeault's theorem.  (See
\S\S\ref{subsec:3.1}--\ref{subsec:3.2}.)  The {\it Euler characteristic\/} is
defined as the alternating sum
  \begin{equation}\label{eq:3}
     \chi (X,V) = \sum\limits_{q=0}^n\, (-1)^q \dim H^q(X,V). 
  \end{equation}
As for the case $n=1$ of Riemann surfaces, one often wants to compute $\dim
H^0(X,V)$, but in general $\dim H^0(X,V)$~ depends on more than topological
data.  On the other hand, the Euler characteristic~$\chi (X,V)$ does have a
topological formula in terms of the Chern classes~$c_j(X)$ and~$c_k(V)$.  The
special case $\dim X = \rank V = 1$ is the classical Riemann-Roch
formula~\eqref{eq:1}.  For $X$~a smooth projective algebraic surface ($n=2$)
and $V\to X$~the trivial bundle of rank~$1$, the result is commonly known as
Noether's formula:
  \begin{equation}\label{eq:4}
     \chi (X) = \frac{1}{12}\bigl( c_1^2(X) + c_2(X)\bigr) [X]. 
  \end{equation}
In~\eqref{eq:4} the Chern classes are evaluated on the fundamental class
of~$X$ given by the natural orientation.  The presence of~12 in the
denominator gives an integrality theorem for the Chern numbers of a
projective surface.
 
The solution to the Riemann-Roch problem for all~$X,V$---that is, the
computation of~\eqref{eq:3}---is one of Hirzebruch's signal achievements.
Hirzebruch's formula is expressed in terms of the Todd polynomials and the
Chern character.  Suppose that the tangent bundle $TX=L_1\oplus \dots \oplus
L_n$ splits as a sum of line bundles, and set $y_i = c_1(L_i) \in
H^2(X;\ZZ)$.  Then the {\it Todd class\/} is
  \begin{equation}\label{eq:5}
     \Todd(X) = \prod_{i=1}^n \frac{y_i}{1-e^{-y_i}}.
  \end{equation}
This is a cohomology class of (mixed) even degree.  Similarly, if $V=K_1\oplus
\dots \oplus K_r$ is a sum of line bundles, with $x_i = c_1(K_i)$, then the
{\it Chern character\/} is  
  \begin{equation}\label{eq:6}
     \ch(V) = \sum\limits_{i=1}^r \,e^{x_i}. 
  \end{equation}
The splitting principle in the theory of characteristic classes allows us to
extend these definitions to ~$TX\to X$ and~$V\to X$ which are not sums of
line bundles.

  \begin{theorem}[Hirzebruch-Riemann-Roch]\label{thm:1}
 Let $X$~be a projective complex manifold and let $V\to X$ be a holomorphic
vector bundle.  Then
  \begin{equation}\label{eq:7}
     \chi (X,V) = \Todd(X)\ch(V)[X]. 
  \end{equation}
  \end{theorem}

Hirzebruch's second main theorem, which is a step in the proof of
Theorem~\ref{thm:1}, is now called Hirzebruch's Signature Theorem.  Let~$X$
be a closed oriented real differentiable manifold of dimension~$4k$ for some
positive integer~$k$.  Then there is a nondegenerate symmetric bilinear
pairing on the middle cohomology~$H^{2k}(X;\RR)$ given by the cup product
followed by evaluation on the fundamental class:
  \begin{equation}\label{eq:8}
     \begin{aligned} H^{2k}(X;\RR) \otimes H^{2k}(X;\RR) &\longrightarrow
     \quad\;\;\,\RR\\ \alpha _1 \quad \;\; \otimes \quad\;\; \alpha_2 \quad\quad
     &\longmapsto 
     (\alpha _1 
     \smallsmile \alpha_2 )[X]\end{aligned} 
  \end{equation}
The signature~$\Sign(X)$ of this pairing is called the signature of~$X$.
(The term `index' is used in place of `signature' in older literature.)
Hirzebruch defines the $L$-class as the polynomial in the Pontrjagin
classes of~$X$ determined by the formal expression
  \begin{equation}\label{eq:9}
     L(X) = \prod_{i=1}^{2k} \frac{y_i}{\tanh y_i}, 
  \end{equation}
where $y_i, -y_i$ are the \emph{Chern roots} of the complexified tangent
bundle.\footnote{The total Pontrjagin class $p(X) = 1 + p_1(X) +
p_2(X)+\cdots $ is defined by the expression $\prod(1 + y_i^2)$.}  This is
analogous to~\eqref{eq:5}: one first defines~$L(X)$ in case $TX\otimes \CC\to
X$~splits as a sum of complex line bundles.

  \begin{theorem}[Hirzebruch Signature Theorem]\label{thm:3}
 The signature of a closed oriented smooth manifold~$X$ is
  \begin{equation}\label{eq:10}
     \Sign(X) = L(X)[X]. 
  \end{equation} 
  \end{theorem}

\noindent
 Hirzebruch's proof uses Thom's bordism theory~\cite{T1} in an essential
way.  Both sides of~\eqref{eq:10} are invariant under oriented bordism and
are multiplicative; for the signature, the former is a theorem of
Thom~\cite[\S IV]{T2}.  Therefore, it suffices to verify~\eqref{eq:10} on a
set of generators of the (rational) oriented bordism ring, which had been
computed by Thom.  The even projective spaces~$\CP^{2n}$ provide a convenient
set of generators, and the proof concludes with the observation that the
$L$-class is characterized as evaluating to~1 on these generators.  The Todd
class enters the proof of Theorem~\ref{thm:1} in a similar manner---its value
on all projective spaces~$\CP^n$ is 1 and it is characterized by this
property.

  \subsection{Grothendieck's Riemann-Roch theorem}\label{subsec:2.2}

The Riemann-Roch-Hirzebruch theorem was extended in a new direction by
Grothendieck~\cite{BS} in 1957.  A decisive step was Grothendieck's
introduction of $K$-theory in algebraic geometry.  Let $X$~be a smooth
algebraic variety.  Define $K(X)$~as the free abelian group generated by
coherent algebraic sheaves on~$X$, modulo the equivalence $\scrF \sim \scrF'
+ \scrF''$ if there is a short exact sequence $0 \to \scrF' \to \scrF \to
\scrF'' \to 0$.  One can replace `coherent algebraic sheaves' by `holomorphic
vector bundles' in this definition, and one fundamental result is that the
group~$K(X)$ is unchanged.  Thus Chern classes and the Chern character are
defined for elements of~$K(X)$.  (Grothendieck refines these to take values
in the {\it Chow ring\/} of~$X$.)  If $f\:X\to Y$ is a morphism of varieties,
and $\scrF$~a sheaf over~$X$, then $R^qf_*(\scrF)$ ~is the sheaf on ~$Y$
associated to the presheaf $U\mapsto H^q(f\inv (U),\scrF)$.  The assignment
  \begin{equation}\label{eq:11}
     f_! \: \scrF\longmapsto \sum (-1)^q R^qf_*(\scrF) \in K(Y) 
  \end{equation}
extends to a homomorphism of abelian groups $f_!\: K(X)\longrightarrow K(Y)$,
as can be seen from the long exact sequence in sheaf cohomology.

Now let $f\:X\to Y$ be a {\it proper\/} morphism between
nonsingular irreducible quasiprojective varieties.  There is a
pushforward~$f_*$ in cohomology (or on the Chow rings).   

  \begin{theorem}[Grothendieck-Riemann-Roch]\label{thm:4}
  For $z\in K(X)$ we have
  \begin{equation}\label{eq:12}
     \Todd(Y)\ch\bigl(f_!(z)\bigr)  = f_*\bigl( \Todd(X)\ch(z)\bigr). 
  \end{equation} 
  \end{theorem}

\noindent
 This reduces to Hirzebruch's Theorem~\ref{thm:1} upon taking $Y$ to be a
point and $z$~the $K$-theory class of a holomorphic vector bundle.

One route to the Todd class is the special case in which $f\:X\to Y$ is the
inclusion of a divisor and $z$~is the class of the structure sheaf~$\sO_X$.
Then $R^qf_*(\sO_X)=0$ for~$q\ge 1$ and $R^0f_*(\sO_X)$~is $\sO_X$~extended
by zero to~$Y$.  Let $L\to Y$~be the line bundle defined by the divisor~$X$.
Observe that $f^*(L)\to X$~is the normal bundle to~$X$ in~$Y$.  The exact
sequence of sheaves
  \begin{equation}\label{eq:13}
     0 \longrightarrow L\inv \longrightarrow \sO_Y \longrightarrow
     f_!\sO_X \longrightarrow 0 
  \end{equation}
leads to the equality
  \begin{equation}\label{eq:14}
     f_!\sO_X = \sO_Y - L\inv 
  \end{equation}
in~$K(Y)$.  
Set $y=c_1(L)$.  Then from~\eqref{eq:14},
  \begin{equation}\label{eq:15}
     \ch\bigl( f_!(\sO_X)\bigr) = 1 - e^{-y}. 
  \end{equation}
On the other hand  
  \begin{equation}\label{eq:16}
     f_*\bigl(\ch(\sO_X)\bigr) = f_*(1) = y. 
  \end{equation}
Thus $f_*\circ \ch = \ch \circ f_!$ up to the Todd class of~$L$.  To check
Theorem~\ref{thm:4} in this case, rewrite~\eqref{eq:12} using the exact
sequence 
  \begin{equation}\label{eq:127}
     0\longrightarrow TX\longrightarrow f^*TY\longrightarrow
     f^*L\longrightarrow 0 
  \end{equation}
of vector bundles on~$X$ and the multiplicativity of the Todd genus: 
  \begin{equation}\label{eq:128}
     \begin{split} \ch\bigl(f_!(z) \bigr) &= f_*\left(
      \frac{\Todd(X)}{f^*\Todd(Y)}\,\ch(z) \right) \\ &= f_*\left(
      \frac{1}{f^*\Todd(L)}\,\ch(z) \right) \\ &=
      \frac{1}{\Todd(L)}\,f_*\bigl(\ch(z) \bigr).\end{split} 
  \end{equation}
This is what we checked in~\eqref{eq:15} and \eqref{eq:16} for~$z=[\sO_X]$.

It is instructive at this stage to consider the inclusion of the zero section
$f\:X\to E$ in a rank~$k$ vector bundle~$\pi \:E\to X$.  Then the
sheaf~$f_!\sO_X=R^0f_*\sO_X $ fits into the exact sequence
  \begin{equation}\label{eq:17}
     0 \longrightarrow \pi ^*{\textstyle\bigwedge} ^kE^* \longrightarrow \pi
     ^*{\textstyle\bigwedge} ^{k-1}E^* \longrightarrow \cdots \longrightarrow \pi
     ^*E^* \longrightarrow \sO_E \longrightarrow f_!\sO_X \longrightarrow 0 
  \end{equation}
of sheaves over~$E$.  (Compare~\eqref{eq:13}.)  Here $E^*$~is the (sheaf of
sections of the) dual bundle to~$\pi \:E\to X$, and the arrows
in~\eqref{eq:17} at~$e\in E$ are contraction by~$e$.  Thus in~$K(E)$ we have 
  \begin{equation}\label{eq:18}
     f_!(\sO_X) = {\textstyle\bigwedge} ^{\bullet }(E^*),
  \end{equation}
where ${\textstyle\bigwedge} ^{\bullet
}(E^*)=\sum\limits_{}(-1)^k{\textstyle\bigwedge} ^kE^*$ in $K$-theory.
Note that $\pi \:E\to X$ is the normal bundle to~$X$ in~$E$.

  \subsection{Integrality theorems in topology}\label{subsec:2.4}

One consequence of the Riemann-Roch-Hirzebruch Theorem~\ref{thm:1} is that
the characteristic number on the right hand side of~\eqref{eq:7}, which {\it
a priori\/} is a rational number, is actually an integer.  This integer is
identified as a sum and difference of dimensions of vector spaces by the left
hand side.  On the other hand, the right hand side is defined for any almost
complex manifold.  Hirzebruch was led to ask (as early as 1954) whether the
{\it Todd genus\/} $\Todd(X)[X]$ of an almost complex manifold (much less a
non-algebraic complex manifold) is an integer~\cite{H3}.  He also asked
analogous questions for real manifolds.  Define the {\it
$\Ahat$-class\/}\footnote{Hirzebruch had previously defined an $A$-class
which differs from the $\Ahat$-class by a power of~2, hence the
notation~$\Ahat$.} of a real manifold~$X^{4k}$ by the formal expression
  \begin{equation}\label{eq:22}
     \Ahat(X) = \prod_{i=1}^{2k} \frac{y_i/2}{\sinh y_i/2},
  \end{equation}
where $y_i$~are the Chern roots.  This is a polynomial in the Pontrjagin
classes.  Then the Todd class of an almost complex manifold can be expressed
as
  \begin{equation}\label{eq:23}
     \Todd(X) = e^{c_1(X)/2} \Ahat(X). 
  \end{equation}
In particular, $\Todd(X)$ depends only on the Pontrjagin classes and the
first Chern class.  It is reasonable to speculate that it was~\eqref{eq:23}
which motivated Hirzebruch to introduce the $\Ahat$-class.  Furthermore,
since the second Stiefel-Whitney class~$w_2$ is the mod~2 reduction of~$c_1$,
Hirzebruch asked: if a real manifold~$X^{4k}$ has $w_2(X)=0$, i.e., if $X$~is
a spin manifold, then is $\Ahat(X)[X]$ an integer?\footnote{In~\cite{H3}
Hirzebruch only asks a less sharp divisibility question (Problem~7 of that
paper).  The more precise form came later, along with the more general
question: If a closed real manifold~$X$ admits an element~$c\in H^2(X;\ZZ)$
whose reduction mod~2 is~$w_2(X)$, and $\Todd(X)$~is defined by~\eqref{eq:23}
(with $c$~replacing $c_1(X)$), then is $\Todd(X)[X]$ an integer?} This was
proved true (initially up to a power of~2 in~\cite{BH2}) by Borel and
Hirzebruch~\cite{BH3} in the late 1950's using results of Milnor on
cobordism~\cite{Mi1}.

The integrality proved, the obvious question presented itself:  
  \begin{equation}\label{eq:24}
     \textnormal{What is the integer~$\Ahat(X)[X]$?} 
  \end{equation}

A first answer to this question came from within algebraic topology, though
not from traditional Eilenberg-MacLane cohomology theory.  When Atiyah and
Hirzebruch learned about Grothendieck's work, they immediately set out to
investigate possible ramifications in topology.  The first step was to define
$K$-theory for arbitrary CW complexes~$X$~\cite{AH1}.  The definition is as
for algebraic varieties, but with `topological vector bundles' replacing
`coherent algebraic sheaves.'  The basic building blocks of topology are the
spheres, and the calculation of~$K(S^n)$ quickly reduces to that of the
stable homotopy groups of the unitary group.  By a fortunate coincidence Bott
had just computed (in 1957) these homotopy groups~\cite{B1}, ~\cite{B2}.  His
{\it periodicity theorem\/} became the cornerstone of the new topological
$K$-theory.  What results is a cohomology theory which satisfies all of the
Eilenberg-MacLane axioms save one, the dimension axiom.  Thus was born
``extraordinary cohomology.''  $K$-theory is the subject of Graeme Segal's
paper in this volume~\cite{Seg4}.

Returning to the Grothendieck program, Atiyah and Hirzebruch formulated a
version of Riemann-Roch for smooth manifolds ~\cite{AH2}, ~\cite{H2}.  Let
$f\:X\to Y$ be a smooth map between differentiable manifolds, and
suppose $f$~is `oriented' in the sense that there exists an element~$c\in
H^2(X;\ZZ)$ with  
  \begin{equation}\label{eq:19}
     c\equiv w_2(X) - f^*w_2(Y) \pmod2. 
  \end{equation}
Recall that Grothendieck's theorem~\eqref{eq:12} is stated in terms of a map
$f_!\:K(X)\to K(Y)$.  In the topological category we cannot push forward
vector bundles, as we could sheaves in the algebraic category, so a new
construction is needed.\footnote{The definition of~$f_!$ was not given in the
original paper~\cite{AH2}; missing was the Thom class in $K$-theory, which is
closely related to the symbol of the Dirac operator.  The Dirac operator
enters the story in the collaboration of Atiyah and
Singer~(\S\ref{subsec:4.1}), and then the $K$-theory Thom class and Thom
isomorphism appear in~\cite[\S12]{ABS}.  See also the discussion
in~\cite[\S1]{Seg4}.} Here we restrict our attention to embeddings of complex
manifolds to simplify the presentation.\footnote{General case: Embed a closed
manifold~ $X$ in a sphere, and so factor an arbitrary map $f\:X \to Y$ into
an embedding followed by a projection: the composition $X \to S^N\times Y \to
Y$.  Bott Periodicity calculates the ``shriek map'' $K(S^N\times Y) \to
K(Y)$.  For embeddings of real manifolds (with an orientation of the normal
bundle) Clifford multiplication on spinors replaces ~\eqref{eq:20}.} Then
\eqref{eq:14} ~ and ~\eqref{eq:18}~ motivate the definition of~$f_!$.  Let
$\pi \:E \to X$ be the normal bundle of~$X$ in~$Y$.  By the tubular
neighborhood theorem, we can identify~$E$ with a neighborhood~$U$ of~$X$
in~$Y$.  The {\it Thom complex\/} ~${\textstyle\bigwedge} ^{\bullet }E^*\to
E$ is defined on the total space of~$E$ by contraction
(compare~\eqref{eq:17}):
  \begin{equation}\label{eq:20}
     0 \longrightarrow \pi ^*{\textstyle\bigwedge} ^kE^*
     \xrightarrow{\;\;\iota (e)\;\;} \pi 
     ^*{\textstyle\bigwedge} ^{k-1}E^* \longrightarrow \cdots \longrightarrow
     \pi ^*E^* \xrightarrow{\;\;\iota (e)\;\;} E\times \CC \longrightarrow 0. 
  \end{equation}
Notice that \eqref{eq:20} is exact for~$e\not= 0$, so the resulting
$K$-theory element is supported on~$X$.  By the tubular neighborhood theorem
it is also defined on~$U$, and extension by zero yields the desired
element~$f_!(1) \in K(Y)$.  If $V \to X$ is a vector bundle, then $f_!(V)$~is
defined by tensoring~\eqref{eq:20} with~$\pi ^*V$.

The Atiyah-Hirzebruch Riemann-Roch theorem for smooth manifolds states
  \begin{equation}\label{eq:21}
     \ch\bigl(f_!(z)\bigr) \Todd(Y) = f_*\bigl( \ch(z)\Todd(X)\bigr), \qquad
     z\in K(X). 
  \end{equation}
Once $f_!$~is defined, the proof is an exercise that compares Thom
isomorphisms in K-theory and cohomology.  Specialize now to $Y=\pt$, and
suppose~$w_2(X)=0$.  Choose the orientation class~$c\in H^2(X;\ZZ)$ to be
zero.  Then for $z=0$ in~\eqref{eq:21} we deduce, in view of~\eqref{eq:23},
that $\Ahat(X)[X] = f_!(1)\in K(\pt)\cong \ZZ$ is an integer.  This argument
by Atiyah-Hirzebruch provided a new proof of the integrality theorem
for~$\Ahat$, and also a topological interpretation of the
integer~$\Ahat(X)[X]$, so a first answer to~\eqref{eq:24}.

Still, that explanation was not considered satisfactory.  As reported by
Atiyah~\cite{A1}, Hirzebruch realized that the signature is the difference in
dimensions of spaces of harmonic differential forms, and he asked for a
similar analytic interpretation of the $\Ahat$-genus $\Ahat(X)[X]$.  Thus
when Singer arrived for a sabbatical stay in Oxford in January, 1962, the
first question Atiyah asked him was, ``Why is the A-roof genus an integer for
a spin manifold?''  Singer~\cite{S1} responded, ``Michael, why are you asking
me that question?  You know the answer to that.''  But Atiyah was looking for
something deeper, and he immediately had Singer hooked.  By March the duo was
in possession of the Dirac operator and the index formula.  Then, nine months
after that initial conversation, Atiyah and Singer completed the first proof
of their eponymous index theorem.

   \section{Antecedents in analysis}\label{sec:3}

The Atiyah-Singer index theorem brings the worlds of algebraic geometry and
algebraic topology together with the worlds of differential geometry and
global analysis.  Our introduction to the latter in~\S\ref{subsec:3.1} begins
with foundational theorems about harmonic differential forms and their
relationship to cohomology.  Geometric elliptic differential operators on
Riemannian manifolds play a central role.  We take up more general elliptic
operators in~\S\ref{subsec:3.2}, where we also recall basic facts about
Fredholm operators.  The Fredholm index, an integer-valued deformation
invariant of a Fredholm operator, is the eponymous character of index theory.
In~\S\ref{subsec:3.3} we give the reader an inkling of the activity around
indices of elliptic operators during the years 1920--1963.

  \subsection{de Rham, Hodge, and Dolbeault}\label{subsec:3.1}

We begin with the de Rham and Hodge theorems, which exemplify the relationship
between elliptic linear differential equations and topology.  Let $X$~be a
smooth $n$-dimensional manifold, and consider the complex of differential
forms
  \begin{equation}\label{eq:25}
     \Omega ^0(X) \xrightarrow{\;\;d\;\;} \Omega ^1(X)
     \xrightarrow{\;\;d\;\;} \cdots \xrightarrow{\;\;d\;\;} \Omega ^n(X), 
  \end{equation}
where $d$~is the exterior derivative of Elie Cartan.  The de Rham cohomology
vector spaces are defined as the quotients
  \begin{equation}\label{eq:26}
     H^q_{\text{DR}}(X) = \frac{\quad\,\ker \,[d\: \Omega ^q(X) \to \Omega
     ^{q+1}(X)]}{\image \,[d\: \Omega ^{q-1}(X) \to \Omega ^q(X)]},\qquad
     0\le q\le n. 
  \end{equation}
The theorem de Rham proved in his 1931 thesis~\cite{deR} states that for
each~$q$ there is a natural isomorphism
$H^q_{\text{DR}}(X)\xrightarrow{\;\;\cong \;\;} H^q(X;\RR)$ of the de Rham
cohomology with real cohomology defined via singular cochains.  (This is
modern language; de Rham proved that there is a closed form with specified
periods, unique modulo exact forms.)  Notice that $H^q_{\text{DR}}(X)$~is
defined using a differential operator, whereas $H^q(X;\RR)$~comes from
topology.  Hodge, motivated by questions in algebraic geometry, proved that
on a {\it closed Riemannian\/} manifold there is a unique ``best'' form in
each cohomology class.  Namely, on an oriented Riemannian manifold~$X$ Hodge
defined a duality operation $*\: \Omega ^q(X)\to \Omega ^{n-q}(X)$, and for
closed manifolds he argued~\cite{Hod} that in each de Rham cohomology class
there is a unique form~$\omega $ satisfying
  \begin{equation}\label{eq:27}
     d\omega =0,\qquad d(*\omega )=0. 
  \end{equation}
These {\it harmonic\/} differential forms are solutions to the elliptic
Hodge-Laplace equation
  \begin{equation}\label{eq:31}
     \Delta \omega = (dd^*+d^*d)\omega =0, 
  \end{equation}
which on a closed manifold is equivalent to the pair of
equations~\eqref{eq:27}.  The number of linearly independent solutions---the
dimension of the vector space~$\sH^q(X)$ of solutions---equals a topological
invariant, the Betti number $b_q=\dim H^q(X;\RR)$.  There is a stronger
statement, namely an isomorphism $\sH^q(X)\xrightarrow{\;\cong
\;}H^q(X;\RR)$.  Neither statement generalizes to arbitrary elliptic
differential operators; rather, the index theorem in this situation computes
the alternating sum of dimensions of spaces of harmonic forms, a familiar
topological invariant:
  \begin{equation}\label{eq:32}
     \sum\limits_{q=0}^n (-1)^q\dim\sH^q(X) = \Euler(X), 
  \end{equation}
where $\Euler(X)$~is the Euler number of~$X$.  (Compare~\eqref{eq:3}.)

We can express the left hand side of~\eqref{eq:32} as the index of an
elliptic operator, namely 
  \begin{equation}\label{eq:77}
     d+d^*\:\Omega ^{\textnormal{even}}(X)\longrightarrow \Omega
     ^{\textnormal{odd}}(X).
  \end{equation}
Its formal adjoint is $d+d^*\:\Omega ^{\textnormal{odd}}(X)\to \Omega
^{\textnormal{even}}(X)$, and we identify the cokernel of~\eqref{eq:77} with
the kernel of the adjoint.  If $\dim X$~is even, then a different
$\zt$-grading\footnote{The grading $\Omega ^{\bullet }(X,\CC)= \Omega
^+(X,\CC)\oplus \Omega ^-(X,\CC)$ is the eigenspace decomposition of the
involution~$\tau $ on~$\Omega ^{\bullet }(X,\CC)$ defined by $\tau (\omega )=
i^{p(p-1)+\ell }*\omega $, where $\dim X=2\ell $.} on complex differential
forms $\Omega ^{\bullet }(X,\CC)$ leads to another elliptic
operator~\cite[\S6]{AS3}
  \begin{equation}\label{eq:78}
     d+d^*\:\Omega ^+(X,\CC)\longrightarrow \Omega ^-(X,\CC), 
  \end{equation}
the \emph{signature operator}, whose index on a closed manifold of dimension
divisible by~4 is the signature of the pairing~\eqref{eq:8}.

Let $X$~be a closed $n$-dimensional complex manifold and $V\to X$ a
holomorphic vector bundle.  Then the sheaf cohomology groups~$H^q(X,V)$ used
in~\S\ref{subsec:2.3} are isomorphic to the cohomology groups of the
$\dbar$-complex
  \begin{equation}\label{eq:33}
     \Omega ^{0,0}(X,V)\xrightarrow{\;\;\dbar\;\;}\Omega
     ^{0,1}(X,V)\xrightarrow{\;\;\dbar\;\;}\cdots \xrightarrow{\;\;\dbar\;\;}
     \Omega ^{0,n}(X,V), 
  \end{equation}
as was proved by Dolbeault.  If $\dim X=1$ then \eqref{eq:33}~reduces to a
single elliptic operator, and for $V\to X$ the line bundle associated to a
divisor~$D$ on~$X$ the vector space~$\sL(D)$ of~\S\ref{subsec:2.1} is
naturally isomorphic to the kernel of~$\dbar$.  If $X$~is K\"ahler, then
Hodge theory implies that the Dolbeault cohomology vector spaces are
isomorphic to vector spaces of complex harmonic differential forms.  Putting
these theorems together, we deduce that on a K\"ahler manifold the
holomorphic Euler characteristic~\eqref{eq:3} is the alternating sum of
dimensions of spaces of harmonic forms.  Hirzebruch's Riemann-Roch
Theorem~\ref{thm:1} is a topological formula for this analytic quantity:
  \begin{equation}\label{eq:35}
     \sum\limits_{q=0}^n \dim\sH^{0,q}(X) = \chi (X,V)=\Todd(X)\ch(V)[X]. 
  \end{equation}
As mentioned at the end of~\S\ref{sec:2}, the Signature Theorem can also be
interpreted in terms of harmonic differential forms.

  \subsection{Elliptic differential operators and the Fredholm
  index}\label{subsec:3.2} 

We set up the index problem on a closed $n$-manifold~$X$.  Let $E^0,E^1\to X$
be vector bundles over ~$X$, and suppose $P\:C^{\infty}(X,E^0) \to
C^{\infty}(X,E^1)$ is a linear differential operator of order~$m$.  In local
coordinates $x^1,\dots ,x^n$ on~$X$, for $u\in C^{\infty}(X,E^0)$ a smooth
section of $E^0\to X$ we have
  \begin{equation}\label{eq:28}
     Pu = a^{i_1i_2\dots i_m} \frac{\partial^mu}{\partial x^{i_1}\partial
     x^{i_2}\cdots \partial x^{i_m}} + \text{lower order terms}, 
  \end{equation}
where $a^{i_1i_2\cdots i_m}$ is a bundle map~$E^0\to E^1$ depending
symmetrically on the~$i_j$, and we sum over the indices~$1\le i_j\le n$.
This highest order term transforms as a tensor under coordinate
changes, so it defines a global bundle map
  \begin{equation}\label{eq:29}
     \sigma (P)\: \Sym^m(T^*X) \otimes E^0\longrightarrow E^1,
  \end{equation}
called the {\it symbol} of~$P$.  View $\sigma (P)$ as a homogeneous
polynomial of degree~$m$ in~$T^*X$ with values in~$\Hom(E^0,E^1)$.  The
differential operator~$P$ is {\it elliptic\/} if its symbol is invertible;
that is, if for each $x\in X$ and nonzero~$\theta \in T^*_xX$, the linear map
$\sigma (P)(\theta ,\dots ,\theta )\:E^0_x\to E^1_x$ ~is invertible.  It
follows from elliptic theory that $P$~has finite dimensional kernel and
cokernel.  (This relies on the \emph{compactness} of~$X$.)  The
(\emph{Fredholm}) \emph{index} of~$P$ is
  \begin{equation}\label{eq:30}
     \ind P = \dim\ker P - \dim\coker P. 
  \end{equation}

Elliptic theory proves that the extension of~$P$ to appropriate Sobolev
spaces is a \emph{Fredholm operator}.  Recall that a Fredholm operator
$H^0\to H^1$ is a bounded linear operator between Hilbert spaces which has
closed range, finite dimensional kernel, and finite dimensional cokernel.
(The definition generalizes to Banach spaces and beyond.)  The index of a
Fredholm operator is defined\footnote{A linear map $T\:V^0\to V^1$ between
finite dimensional vector spaces is an element of $(V^1)\otimes (V^0)^*$, so
it stands to reason that the sign in~\eqref{eq:30} should have been the
opposite: the dual---or minus---is the domain, not the codomain.  The usual
sign convention causes headaches down the road, for example in the theory of
determinants.  On the other hand, one could argue for the usual sign
convention by rewriting a single operator as a 2-term complex in degrees~0
and~1, and then the usual sign for the index (as in~\eqref{eq:3}) reduce to~
\eqref{eq:30}.}  by~\eqref{eq:30}.  The space $\Hom(H^0,H^1)$ of continuous
linear maps has a Banach space structure defined by the operator norm, and
the open subspace $\Fred(H^0,H^1)\subset \Hom(H^0,H^1)$ of Fredholm operators
has nontrivial homotopy groups of unbounded degree.  In particular, the index
  \begin{equation}\label{eq:34}
     \ind \:\pi _0\Fred(H^0,H^1)\longrightarrow \ZZ 
  \end{equation}
is an isomorphism.  In other words, the numerical index is a complete
deformation invariant of a single Fredholm operator.  Furthermore, the action
of compact operators by translation preserves the subspace of Fredholm
operators, hence the index is invariant under this translation.

Elliptic theory implies that the lower order terms of an elliptic
differential operator~\eqref{eq:28} on a smooth manifold are compact relative
to the highest order term, which is essentially the symbol~\eqref{eq:29}.  It
follows that the index of an elliptic differential operator is an invariant
of its symbol.  Furthermore, a continuous path of elliptic differential
operators induces a continuous path of Fredholm operators and of symbols (in
suitable topologies).  The index is unchanged under such deformations and, in
fact, only depends on the homotopy class of the symbol.  The Atiyah-Singer
index theorem provides a formula for the index in terms of the homotopy class
of the symbol.

  \subsection{Index problems for elliptic operators}\label{subsec:3.3}

There is a long and rich history of index theorems for linear elliptic
problems in the first half of the $20^{\textnormal{th}}$~century.  Many are
subsumed by the Atiyah-Singer index theorem and its extension to manifolds
with boundary~(\S\ref{subsec:5.2}).  The articles by Agranovich~\cite{Agr}
and Seeley~\cite{Se1} are excellent guides to this history.  The first index
theorem is contained in a 1920~paper of Fritz Noether~\cite{N}.  (This is
credited in modern references such as~\cite{E,AM}.)  Moreover, this paper
seems to be the origin of the Fredholm index.  In fact, in the older
literature the following terminology that is sometimes used: a linear
operator is said to obey the `Noether property' if it is Fredholm, in which
case `Fredholm' is reserved for an operator of index zero (which then
satisfies the ``Fredholm alternative'').
 
One case of Noether's work is an index formula for Toeplitz operators.  Let
$\cir\subset \CC$ be the unit circle and $f\:\cir\to \CC^\times $ a smooth
nonzero complex-valued function.  The Toeplitz operator $T_f\:\sH\to \sH$ is
defined on the Hilbert space~$\sH$ of $L^2$~holomorphic functions on the
closure~$\overline{\Omega }$ of the unit disk~$\Omega \subset \CC$.  By
Fourier series $\sH$~sits as a subspace in $\tsH=L^2(\cir,\CC)$.  Let
$i\:\sH\hookrightarrow \tsH$ be the inclusion and $\pi
\:\tsH\twoheadrightarrow\sH$ the orthogonal projection.  Then $T_f=\pi \circ
M_f\circ i$ is the compression of the multiplication operator $M_f\:\tsH\to
\tsH$ to~$\sH$.   

  \begin{theorem}[Noether, 1920]\label{thm:32}
 $T_f$~is Fredholm with index minus the winding number of~$f$. 
  \end{theorem}

\noindent
 The reader may wish to compute the index explicitly for $f_n(z)=z^n$, $n\in
\ZZ$.  Theorem~\ref{thm:32}~was rediscovered by Gohberg-Krein~\cite{GK}, and
there are index theorems for Toeplitz operators in arbitrary dimensions, for
example in Boutet de Monvel~\cite{BdM}.

In 1960, Gelfand~\cite{G} observed that the index is a homotopy invariant,
and he posed the general problem of computing a topological formula for the
index.  It seems that Atiyah and Singer were unaware of these events in
Russia when they embarked on the journey which led to the index theorem,
though they became aware of them during a visit by Smale to Oxford~\cite{A1}.
Gelfand's paper, and some of its antecedents which solve special cases of the
index problem, are cited at the beginning of the Atiyah-Singer announcement
of their general index theorem~\cite{AS1}.

   \section{The index theorem and proofs}\label{sec:4}

We arrive at the Atiyah-Singer index theorem for a single elliptic operator.
It was announced in~\cite{AS1} in~1963.  Atiyah-Singer's first proof
(\S\ref{subsec:4.2}), modeled on Hirzebruch's cobordism proofs of his
signature and Riemann-Roch theorems, was not written up by them but rather
was published in Palais~\cite{Pa} as a series of pieces by many contributors
in a volume which remains a valuable reference.  The second proof
(\S\ref{subsec:4.4}), modeled more on Grothendieck, appeared in ~1968 in the
first~\cite{AS2} of a series of papers by Atiyah and Singer.  Subsequent
papers treat variations and generalizations.  We begin in~\S\ref{subsec:4.1}
with the Dirac operator in Riemannian geometry.  It is the analogue of Dirac's
operator in Lorentz geometry, and it is central in many contexts in geometry
and physics as well as in general index theory.  In a different direction,
pseudodifferential operators play an important role in both proofs of the
index theorem; we give a brief introduction in~\S\ref{subsec:4.3}.
In~\S\ref{subsec:4.5} we list a few early applications of the index theorem.

  \subsection{The Dirac operator}\label{subsec:4.1}

In~1928, Dirac~\cite{D} introduced his equation as part of his relativistic
theory of electrons.  Dirac worked on Minkowski spacetime.  The analogue of
Dirac's line of inquiry for Euclidean space~$\EE^n$ with standard coordinates
$x^1,\dots ,x^n$  asks for a first-order differential operator
  \begin{equation}\label{eq:36}
     D = \gamma ^1\frac{\partial }{\partial x^1} +\cdots+\gamma
     ^n\frac{\partial }{\partial x^n} 
  \end{equation}
whose square is the Laplace operator 
  \begin{equation}\label{eq:37}
     \Delta =-\left\{ \left(\frac{\partial }{\partial x^1}\right)^2 + \cdots
     + \left(\frac{\partial }{\partial x^n}\right)^2 \right\} . 
  \end{equation}
(In Minkowski spacetime the elliptic Laplace operator~\eqref{eq:37} is
replaced by the hyperbolic wave operator.)  Assume that $\gamma ^1,\dots
,\gamma ^n$~are constant functions on~$\EE^n$.  Then the differential
equation
  \begin{equation}\label{eq:38}
     D^2=\Delta 
  \end{equation}
is equivalent to the system of algebraic equations 
  \begin{equation}\label{eq:39}
     \gamma ^i\gamma ^j + \gamma ^j\gamma ^i = -2\delta ^{ij} = \begin{cases}
     -2,&i=j;\\\phantom{-}0,&i\neq j,\end{cases}\qquad 1\le i,j\le n. 
  \end{equation}
There are no scalar solutions if~$n\ge2$, but there are matrix solutions.
Let\footnote{To define the closely related Clifford algebra~$\Cnp$, change
the $-$~ sign in relation~\eqref{eq:39} to a $+$~sign.  Both Clifford
algebras appear in~\S\ref{subsec:5.6}.} $\Cnm$, the \emph{Clifford
algebra}~\cite{Cl,Ca,BW,Ch,ABS}, be the unital algebra over~$\RR$ generated
by $\gamma ^1,\dots ,\gamma ^n$ subject to the relations~\eqref{eq:39}.  A
matrix solution to~\eqref{eq:39} defines a $\Cnm$-module.  The \emph{spin
group}~$\Spin_n$ is a double cover of the special orthogonal group~$\SO_n$,
and $\Spin_n$~ is a subgroup of the units in~ $\Cnm$, much as $\SO_n$~is a
subgroup of the units in the algebra of $n\times n$~real matrices.  Clifford
modules restrict to special representations of~$\Spin_n$ called \emph{spin}
or \emph{spinor} representations.
 
Why does the spin group enter the quest to identify the integer
$\Ahat(X)[X]$?  One key is the formula~\eqref{eq:22} for the $\Ahat$-genus.
Assuming $n$~is even, the group~ $\Spin_n$ has two distinguished inequivalent
irreducible complex representations~$\SS^0,\SS^1$, and the difference of the
characters of~$\SS^0$ and~$\SS^1$ is an upside down variant of the
$\Ahat$-genus~\eqref{eq:22}, namely
  \begin{equation}\label{eq:40}
     \prod\limits_{i=1}^{n/2}\frac{\sinh y_i/2}{1/2}, 
  \end{equation}
where now $y_1,\dots ,y_n$ are a basis of characters of a maximal torus
of~$\Spin_n$.  Second, the spin condition was known to be related to the
integrality of the $\Ahat$-genus, as explained in~\S\ref{subsec:2.4}.
Perhaps these considerations led Atiyah and Singer to construct the Dirac
operator on a Riemannian spin manifold.  Observe too that on a Riemannian
manifold the square of the first-order differential operator $d+d^*$ is the
Hodge-Laplace operator~\eqref{eq:31}, so in that sense $d+d^*$ is already a
Dirac operator.  On a K\"ahler manifold, the same holds for $\dbar + \dbar^*$.
 
\emph{The} Dirac operator on an even dimensional Riemannian spin
manifold~$X^n$ is a main character in the index theorem, so we give the
definition here.  Let $\SO(X)\to X$ be the principal $\SO_n$-bundle whose
fiber at~$x\in X$ consists of oriented isometries $\RR^n\to T_xX$.  The spin
structure is a double cover $\Spin(X)\to\SO(X)$ together with a compatible
principal $\Spin_n$-bundle structure on the composition $\Spin(X)\to\SO(X)\to
X$.  Use the complex spin representations~$\SS^0,\SS^1$ to construct
associated complex vector bundles $\SS^{0}_X,\SS^{1}_X\to X$; sections of
these vector bundles are called \emph{spinors} or \emph{spinor fields}.  The
Levi-Civita connection on~$X$ induces a covariant derivative~$\nabla $ on
each spinor bundle.  The Dirac operator is
  \begin{equation}\label{eq:41}
     D_X = c\circ \nabla \:C^{\infty}(X,\SS^0_X)\longrightarrow
     C^{\infty}(X,\SS^1_X),  
  \end{equation}
where $c$~is Clifford multiplication, induced from a $\Spin_n$-equivariant
linear map $\RR^n\otimes \SS^{0}\to\SS^{1}$.  There is also a Clifford
multiplication $\RR^n\otimes \SS^1\to \SS^0$, and both together give
$\SS^0\oplus \SS^1$ the structure of a $\Cliff_n$-module.  This construction
is not quite canonical, since the irreducible representations~$\SS^0,\SS^1$
are only determined uniquely up to tensoring by a line.
In~\S\ref{subsec:5.6} we introduce the Clifford linear Dirac operator, which
is canonical.
 
The de Rham~\eqref{eq:25} and Dolbeault~\eqref{eq:33} operators have similar
descriptions, and so \eqref{eq:32}~ and \eqref{eq:35} ~ motivate a
conjectural index formula for the Dirac operator.  Let $\sH^{\pm}(X)$ be the
complex vector spaces of \emph{harmonic spinors}, i.e., solutions~$\psi $ to
$D\psi =0$.  The conjectured index formula is
  \begin{equation}\label{eq:42}
     \dim\sH^+(X) - \dim\sH^-(X) = \Ahat(X)[X]. 
  \end{equation}
For an investigation into various aspects of harmonic spinors, see the DPhil
thesis~\cite{Hi} of Atiyah's student Nigel Hitchin.

  \subsection{First proof: cobordism}\label{subsec:4.2}
 
Let $X$~be a closed oriented manifold; $E^0,E^1\to X$ complex vector bundles;
and $P\:C^{\infty}(X,E^0)\to C^{\infty}(X,E^1)$ an elliptic differential
operator.  The \emph{analytic index} of~$P$ is the Fredholm index defined
in~\eqref{eq:30}.  The \emph{topological index}---in cohomological form---is
the following.  Ellipticity implies that the symbol~$\sigma (P)$
in~\eqref{eq:29} restricted to the nonzero vectors $T^*X\setminus 0$ is an
isomorphism.  Therefore, $\sigma (P)$~defines a relative $K$-theory class in
$K^0(T^*X,\,T^*X\setminus 0)$.  The Chern character maps $K$-theory to
rational cohomology, and then the inverse of the Thom isomorphism
  \begin{equation}\label{eq:43}
     \phi \:H^{\bullet }(X;\QQ)\xrightarrow{\;\;\cong \;\;} H^{\bullet +\dim
     X}(T^*X,\,T^*X\setminus 0;\,\QQ) 
  \end{equation}
brings us to the cohomology of the base~$X$.  The \emph{topological index} is
the cup product of the Todd class of the complexified tangent bundle with the
image of~$\sigma (P)$ under the Chern character and Thom isomorphism.  The
index theorem asserts that the analytic and topological indices are equal.

  \begin{theorem}[Atiyah-Singer Index Theorem~\cite{AS1}]\label{thm:5}
 The index of~$P$ is 
  \begin{equation}\label{eq:44}
     \ind P = \bigl(\Todd(TX\otimes \CC)\smallsmile \phi \inv \ch \sigma
     (P)\bigr)[X].  
  \end{equation} 
  \end{theorem}

\noindent
 In the remainder of this section we sketch the main ideas which enter the
proof of Theorem~\ref{thm:5}.

  \begin{remark}[]\label{thm:6}
 The $K$-theoretic formula for the index (\S\ref{subsec:4.4}) is more natural
and lends itself to many generalizations.  This fits with Atiyah's philosophy
that $K$-theory, based on linear algebra, is more elementary than cohomology.
Certainly it is the form of algebraic topology which most closely matches
linear differential operators.
  \end{remark}

Now to the proof.  At the end of~\S\ref{subsec:3.2} we indicated that the
analytic index depends only on the homotopy class of the symbol~$\sigma (P)$.
Atiyah-Singer introduce the abelian group
  \begin{equation}\label{eq:45}
     \Ell(X) = K^0(T^*X,\,T^*X\setminus 0) 
  \end{equation}
of elliptic symbol classes.  It is a module over the ring~$K^0(X)$, by tensor
product, and the \emph{cyclic} module $K^0(X)\sigma _0$ generated by the symbol
class of the signature operator is a subgroup of finite index.  This reduces
the problem on a fixed manifold~$X$ to the signature operator twisted by a
vector bundle $W\to X$.  However, to carry through this argument one needs
that every element of~$\Ell(X)$ is the symbol of an elliptic operator, and
furthermore if the symbols of elliptic operators~$P_0,P_1$ define the same
element of~$\Ell(X)$, then there is a homotopy~$P_0\sim P_1$.
\emph{Differential} operators are too rigid for these properties to hold, and
a critical move in the proof is the introduction of \emph{pseudodifferential}
operators, which we discuss briefly in the next section.

  \begin{remark}[]\label{thm:7}
 If $X$~is a spin manifold, then $\Ell(X)$ is the cyclic $K^0(X)$-module
generated by the symbol of the Dirac operator.  This fact is an expression of
Bott periodicity, as realized in $K$-theory.
  \end{remark}

Bordism enters the proof at this stage.  For the signature operator twisted
by a vector bundle $W\to X$, both sides of~\eqref{eq:44} are viewed as
functions of a pair~$(X,W)$, where we only use the equivalence class $[W]\in
K^0(X)$.  A crucial step is the proof that each side is a bordism invariant
of~$(X,W)$.  This is straightforward for the cohomological formula on the
right hand side.  For the analytic index, suppose $(X,W)$ is the boundary of
a pair~$(Y,U)$, where $Y$~is a compact oriented manifold and $U\to Y$ a
complex vector bundle.  Atiyah-Singer introduce an elliptic differential
operator~$D_U$ on~$Y$ whose boundary operator on~$X$ is the twisted signature
operator.  They specify a \emph{local} elliptic boundary condition, and prove
(i)~$\ind_YD_U = \ind_X D_W$ and (ii)~$\ind_YD_U=0$.  With  bordism invariance
in hand, it remains to compute a basis for the rational vector space of
bordism classes of pairs~$(X,W)$ and verify~\eqref{eq:44} for those.

The journey from this rough sketch to a complete proof is replete with
interesting detours in analysis, geometry, and topology.

  \subsection{Pseudodifferential operators}\label{subsec:4.3}

Let $P$~be a differential operator of order~$m$ on~$\RR^n$, as
in~\eqref{eq:28}.  Its action on a smooth function~$u$ of compact support is
conveniently written in terms of the Fourier transform~$\hat{u}$:
  \begin{equation}\label{eq:46}
     Pu(x) = (2\pi )^{-n}\int_{{\RR^n}^*} p(x,\xi )\,\hat{u}(\xi )\,e^{i\langle
     x,\xi \rangle} \,d\xi ,\qquad x\in \RR^n, 
  \end{equation}
where for each~$x\in \RR^n$ the function $p(x,\xi )$ is a polynomial of
degree~$m$ in~$\xi \in {\RR^n}^*$.  A generalization, going back to
Mikhlin~\cite{Mik} and Calder\'on-Zygmund~\cite{CZ}, allows more general
(total) symbols~$p$, with uniform bounds on the behavior of~$p$ and its
derivatives as $|\xi |\to\infty $.  One motivation comes from elliptic
operators: a parametrix---an approximate inverse---is a linear operator of
this form.  Also, these operators have a Schwartz kernel which is smooth away
from the diagonal, and this leads to good properties.  (The Schwartz kernel
of a differential operator is zero off the diagonal.)  The theory of
pseudodifferential operators, and their globalization to smooth manifolds, is
treated in papers of Seeley~\cite{Se2,Se3}, H\"ormander~\cite{Ho2,Ho3},
Kohn-Nirenberg~\cite{KN}, and Palais-Seeley~\cite{PS} among others.  This is
only a very small sample of the extensive literature.
 
If $P$~is a differential operator of order~$m$, so $p(x,\xi )$~is a
degree~$m$ polynomial in~$\xi $ for each~$x$, then the \emph{principal} or
\emph{top order} symbol of~$P$ is the homogeneous polynomial of degree~$m$
given as 
  \begin{equation}\label{eq:47}
     \sigma (P)(x,\xi ) = \lim\limits_{\lambda \to\infty }\frac{p(x,\lambda
     \xi )}{\lambda ^m}.
  \end{equation}
Restrict to pseudodifferential operators with symbol~$p$ for which the
limit~\eqref{eq:47} exists.  Then the principal symbol is defined, in global
form~\eqref{eq:29} on a smooth manifold, and ellipticity is as before
invertibility of the principal symbol.  The surjectivity and continuity of
the principal symbol map on elliptic pseudodifferential operators are crucial
ingredients in the proof of the index theorem.
 
We remark that in the Atiyah-Bott~\cite{AB1} and Atiyah-Singer~\cite{AS2,AS4}
papers on index theory, the global theory of pseudodifferential operators is
expanded further.

  \subsection{A few applications}\label{subsec:4.5}

One immediate consequence of the Atiyah-Singer Index Theorem~\ref{thm:5} is
the index formula for the Dirac operator~\eqref{eq:42}.  This provides an
\emph{analytic} interpretation of the $\Ahat$-genus of a spin manifold, hence
an answer to~\eqref{eq:24}.
 
Several additional consequences are described in~\cite{AS1}.  As already
mentioned, Hirzebruch's Signature Theorem~\ref{thm:3} is a special case.  So
is the Hirzebruch-Riemann-Roch Theorem~\ref{thm:1}, but the more flexible
techniques of Atiyah-Singer prove it for arbitrary closed complex manifolds,
a powerful generalization from projective algebraic manifolds.  Finally, for
systems of elliptic operators on~$X$---for $E^0,E^1\to X$ trivial vector
bundles---the index vanishes given appropriate inequalities between~$\dim X$
and $\rank E^0=\rank E^1$, a result which connects Theorem~\ref{thm:5} to the
PDE literature of the period.  In particular, the index vanishes for an
elliptic operator acting on a single function.

As new index theorems proliferate, so too do applications, as we will see
in~\S\ref{sec:5}.

  \subsection{Second proof: $K$-theory}\label{subsec:4.4}

Recall the pushforward~\eqref{eq:11} which occurs in Grothendieck's version
of the Riemann-Roch theorem.  Let $X$~be a compact projective variety.  For
the special case of the unique map $f\:X\to\pt$, the pushforward
$f_!\:K(X)\to K(\pt)\cong \ZZ$ computes the Euler characteristic~\eqref{eq:3}
of a holomorphic vector bundle $V\to X$.  According to the Dolbeault Theorem,
the sheaf cohomology groups are isomorphic to the cohomology groups of the
$\dbar$-complex~\eqref{eq:33}.  For smooth manifolds and smooth maps,
Atiyah-Hirzebruch found a \emph{topological pushforward} in one of their
first works on $K$-theory, and it is designed to match Grothendieck's~$f_!$
in this situation.  (See the text surrounding~\eqref{eq:19}.)  There is also
an analytic interpretation.  Namely, on a K\"ahler manifold~$X$ the Dolbeault
cohomology group~$H^{0,q}_{\dbar}(X,V)$ is isomorphic to the kernel of the
elliptic operator $\dbar + \dbar^*$ on $(0,q)$-forms, hence the Euler
characteristic~\eqref{eq:3} is the index of
  \begin{equation}\label{eq:48}
     \dbar+\dbar^*\:\Omega ^{0,\textnormal{even}}(X,V) \longrightarrow \Omega
     ^{0,\textnormal{odd}}(X,V). 
  \end{equation}
That index is an \emph{analytic pushforward} of $V\to X$ under the map~$f$.
The Riemann-Roch-de Rham-Hodge-Dolbeault-Hirzebruch theorems imply the
equality of analytic and topological pushforwards.

The $K$-theory form of the Atiyah-Singer Index Theorem is a generalization
for arbitrary elliptic (pseudo)differential operators.  Recall
from~\S\ref{subsec:4.2} that the ``homotopy class'' of the symbol~$\sigma
(P)$ of an elliptic differential operator~$P$ on a smooth manifold~$X$ is an
element $[\sigma (P)]\in K^0(T^*X, T^*X\setminus 0)$.  Atiyah-Singer define two
homomorphisms
  \begin{equation}\label{eq:49}
     \ind\:K^0(T^*X, T^*X\setminus 0)\longrightarrow \ZZ. 
  \end{equation}
The \emph{analytic index} a-ind takes a symbol class~$\sigma $ to the
index~$\ind P$ of any elliptic pseudodifferential operator~$P$ with $[\sigma
(P)]=\sigma $.  The \emph{topological index} t-ind is similar to the
topological pushforward.  It is based on the Thom isomorphism in $K$-theory,
which in turn rests on Bott periodicity.

  \begin{theorem}[Atiyah-Singer~1967]\label{thm:8} a-ind = t-ind. 
  \end{theorem}

\noindent
 This is equivalent to the cohomological Theorem~\ref{thm:5}, but as we will
see in~\S\ref{sec:5} the naturalness of the $K$-theory formulation and proof
lend themselves to many generalizations. 
 
Theorem~\ref{thm:8} is the subject of~\cite{AS2}.  There is a concise summary
of the proof idea in~\S1 of that paper, though the actual proof follows a
slightly different arc.  In essence, Atiyah-Singer uniquely characterize
``index homomorphisms''~\eqref{eq:49} by a short list of axioms, which they
verify a-ind and t-ind satisfy.  Beyond normalization axioms, two main
properties of the index feature in the proof.  An \emph{excision axiom}
extends the index to compactly supported symbols on arbitrary (potentially
noncompact) manifolds~$U$, at least if $U$~is embeddable in a compact
manifold, by asserting the independence of the index of the extension by zero
of a compactly supported symbol on~$U$ to a symbol on an ambient compact
manifold.  A \emph{multiplicative axiom} tells a product formula for
(twisted) product symbols.  A robust global theory of pseudodifferential
operators is used at this point in the proof.

   \section{Variations on the theme}\label{sec:5}

After the initial Atiyah-Singer work, index theory branched out in multiple
directions.  Atiyah was at the center of many of these developments, which we
can only touch upon in this section.  We begin in~\S\ref{subsec:5.1} with
Atiyah-Bott's generalization of the classical Lefschetz fixed point theorem
to elliptic complexes.  Its manifold applications include a geometric proof
of Weyl's character formula in the theory of compact Lie groups.  The link
between the Atiyah-Singer $K$-theoretic framework for index theory and the
index theory for elliptic operators which preceded it is the extension of the
Atiyah-Singer theorem to compact manifolds with boundary, carried out jointly
with Bott as we recount in~\S\ref{subsec:5.2}.  The basic Atiyah-Singer
theorem applies to complex elliptic operators.  There is an important
refinement~(\S\ref{subsec:5.3}) to {real} elliptic operators which,
naturally, links to the real version of $K$-theory introduced by Atiyah at
this time.  There is a \emph{canonical} real Dirac operator on a Riemannian
spin manifold---perhaps not as oft-used as it deserves to be---so we take
the opportunity in~\S\ref{subsec:5.6} to expose it.  Grothendieck's dictum to
do geometry over a base, not just over a point, leads to an index theorem
for families of elliptic operators.  As Grothendieck's philosophy promises,
this theorem, described in~\S\ref{subsec:5.4}, appears often in geometry and
physics.  Finally, in~\S\ref{subsec:5.5} we introduce an index theorem of
Atiyah's which launched an entirely new branch of the subject, though not one
he himself pursued.

  \subsection{Equivariant index theorems}\label{subsec:5.1}

Recall the classical Lefschetz fixed point theorem.  Let $X$~be a smooth
compact manifold and $f\:X\to X$ a diffeomorphism with isolated fixed points.
The subset $\Fix(f)\subset X$ is then finite, and our aim might be to compute
its cardinality.  But the cardinality is not a deformation invariant---under
deformation pairs of fixed points undergo birth and death---hence there is
not a topological formula.  However, if we count with signs, weighting each
$x\in \Fix(f)$ by\footnote{under the assumption that $df_x$~has no nonzero
fixed vectors for all~$x\in \Fix(f)$}
  \begin{equation}\label{eq:50}
     \epsilon _x=\sign\det(1-df_x)\,\in \{-1,+1\}, 
  \end{equation}
then the Lefschetz fixed point theorem asserts 
  \begin{equation}\label{eq:51}
     \sum\limits_{x\in \Fix(f)}\epsilon _X = \sum\limits_{q=0}^{\dim X}
     (-1)^q\Tr\Bigl( f^*\res{H^q(X;\RR)} \Bigr) , 
  \end{equation}
where the right hand side is the alternating sum of traces
of~$f$ acting on the $q^{\textnormal{th}}$~real cohomology.

In~1965, Atiyah and Bott formulated and proved a generalization to elliptic
operators and elliptic complexes.  Let $X$~be a closed manifold; $E^0,E^1\to
X$ vector bundles; $f\:X\to X$ a diffeomorphism with isolated fixed points;
$f^i\:E^i\to E^i$, $i=0,1$, lifts of~$f$; and $P\:C^{\infty}(X,E^0)\to
C^{\infty}(X,E^1)$ an elliptic differential operator which commutes
with~$f^i$.  Define 
  \begin{equation}\label{eq:52}
     \nu _x = \frac{\Tr\Bigl( f^0\res{E^0_x} \Bigr) - \Tr\Bigl(
     f^1\res{E^1_x} \Bigr) }{\bigl| \det(1-df_x) \bigr| },\qquad x\in
     \Fix(f). 
  \end{equation}

  \begin{theorem}[Atiyah-Bott~\cite{AB1}]\label{thm:9}
 We have 
  \begin{equation}\label{eq:53}
     \sum\limits_{x\in \Fix(f)}\nu _x = \Tr\Bigl(f^0\res{\ker P}\Bigr) -
     \Tr\Bigl( f^1\res{\coker P}\Bigr). 
  \end{equation}
  \end{theorem}

\noindent 
 We will not comment on the proof here, which also involves
pseudodifferential operators, but instead tell some applications.
 
In their proof of the index theorem modeled after
Grothendieck~(\S\ref{subsec:4.5}), Atiyah-Singer~\cite{AS2} incorporate the
action of a \emph{compact} Lie group~$G$ on an elliptic operator~$P$ on a
compact manifold~$X$.  The Atiyah-Bott setup overlaps in case their
diffeomorphism~$f$ of~$X$ generates a (compact) torus of diffeomorphisms.  In
the Atiyah-Singer case, since $G$~commutes with~$P$, it acts on~$\ker P$ and
$\coker P$, so the formal difference---the index---is a virtual
representation of~$G$.  Hence it defines a class in the representation
ring~$R(G)$.  On the other hand, the symbol~$\sigma (P)$ also commutes
with~$G$, and so it defines a class in \emph{equivariant $K$-theory}, an
important extension of topological $K$-theory which was investigated by
Atiyah's student Graeme Segal in his DPhil thesis~\cite{Seg1}.  The topological
index t-ind, executed in equivariant $K$-theory, also produces an element
of~$R(G)$.  The equivariant extension of Theorem~\ref{thm:8} identifies these
elements of~$R(G)$.  Atiyah-Segal~\cite{ASeg} apply the localization theorem
in equivariant $K$-theory~\cite{Seg1} to derive a Lefschetz-type formula, and
in particular they recover Theorem~\ref{thm:9} in case $f$~generates a
compact Lie group (torus) of diffeomorphisms.
 
These various Lefschetz formul\ae\ have many applications; those
contemporaneously realized with the theorems are presented in~\cite{AB2}
and~\cite{AS3}.  Here is a small sample.  Let $X$~be a connected closed
complex manifold with $H^q(X;\sO_X)=0$ for~$q>0$; then any holomorphic map
$f\:X\to X$ has a fixed point.  In a different direction, let $G$~be a
compact connected Lie group and $T\subset G$ a maximal torus.  Then
Theorem~\ref{thm:9} applied to the action of a generic element of~$T$ on a
holomorphic line bundle $\mathscr{L}\to G/T$ over the flag manifold leads to
the Weyl character formula.  On a compact Riemannian manifold~$X$ one deduces
strong consequences of the generalized Lefschetz theorems using the signature
operator.  For example, Atiyah-Bott-Milnor prove that two $h$-cobordant lens
spaces are isometric.  Some number theoretic aspects of the Lefschetz formula
for the signature operator are described in the Atiyah-Bott and Atiyah-Singer
papers.  These ideas are elaborated much further by
Hirzebruch-Zagier~\cite{HZ}, who explore this Lefschetz signature theorem not
only for lens spaces, but also for projective spaces, Brieskorn varieties,
and algebraic surfaces.  They find deep relations with cotangent sums, the
Dedekind eta function, modular forms, and real quadratic fields, among other
number theoretic objects of interest.  This link between Lefschetz-type
invariants and number theory is central to the developments in the 1970s, as
we take up in~\S\ref{sec:7}.

  \subsection{Index theorem on manifolds with boundary}\label{subsec:5.2}

As discussed in~\S\ref{subsec:3.3}, the study of linear elliptic equations
traditionally takes place on a domain~$\Omega \subset \AA^n$ with smooth
boundary, in which case \emph{elliptic boundary conditions} need be imposed.
($\AA^n$~is the standard real affine $n$-dimensional space.) For the
second-order Laplace operator~\eqref{eq:37} on~$\overline{\Omega }$, the
\emph{Dirichlet problem} is the system of equations
  \begin{equation}\label{eq:54}
     \begin{aligned} \Delta u&=0 \\ u\res{\partial \Omega }&=f\end{aligned} 
  \end{equation}
for a prescribed function~$f$ on~$\partial \Omega $.  The second equation
in~\eqref{eq:54} is a \emph{local} boundary condition for~$\Delta $: at each
point~$x\in \partial \Omega $ it depends only on~$u(x)$.  More generally, a
local boundary condition may depend on a finite set of derivatives of~$u$
at~$x$.  Lopatinski~\cite{Lo} gave a general criterion for a local boundary
condition to be elliptic; see~\cite[\S20.1]{Ho1}, for example. 
 
Not all elliptic operators admit local elliptic boundary conditions.
Consider the first-order $\dbar$-operator on functions $u\:\overline{\Omega
}\to\CC$ where $\Omega \subset \CC$ is the unit disk.  The kernel of~$\dbar$
consists of holomorphic functions on~$\overline\Omega $, an infinite
dimensional vector space, as expected in the absence of elliptic boundary
conditions.  For this operator there do not exist \emph{local} elliptic
boundary conditions: the Cauchy-Riemann equations on a disk are ill-posed.
(As we will see in~\S\ref{subsec:7.2}, there are \emph{global} boundary
conditions.  Also, if we consider the $\dbar$-operator on sections of a
vector bundle---as, for example, when studying deformations of holomorphic
disks in a symplectic manifold with boundary on a Lagrangian
submanifold---then there do exist local boundary conditions.)  The
topological nature of the existence question for local elliptic boundary
conditions was brought out in work of Atiyah-Bott~\cite{AB3,A2}.  Again,
$K$-theory is the natural home for the obstruction.  In fact, their work led
to a new proof of the Bott Periodicity Theorem~\cite{AB4}.  Let $X$~be a
compact manifold with boundary, and suppose $\sigma \in
K^0(T^*X,\,T^*X\setminus 0)$ is an elliptic symbol class.  The obstruction to
local elliptic boundary conditions is the restriction of~$\sigma $
to~$K^0(T^*X\res{\partial X})$, and if the obstruction vanishes then local
elliptic boundary conditions lift~$\sigma $ to an appropriate relative
$K$-theory group.  Atiyah-Bott-Singer prove an index theorem for elliptic
operators with local elliptic boundary conditions.  The topological index
easily extends to the lifted relative elliptic symbol, and the topological
index of that lifted symbol computes the analytic index of the elliptic
boundary value problem.

  \subsection{Real elliptic operators}\label{subsec:5.3}

Consider the operator $D=d/dx$ acting on real-valued functions on~$\RR/2\pi
\ZZ$, i.e., on $2\pi $-periodic functions $u\:\RR\to\RR$.  Integration by
parts shows that $D$~is (formally) skew-adjoint.  It is an example of a
\emph{real} Dirac operator.  A closely related Dirac operator is $D'=d/dx$
acting on twisted functions on the circle, or equivalently $u\:\RR\to\RR$
which satisfy $u(x+2\pi )=-u(x)$ for all~$x\in \RR$.  Observe $\ker
D\cong \RR$ consists of constant functions, whereas $\ker D'=0$.  A
skew-adjoint operator has integer index zero in the sense of~\eqref{eq:30}.
But for skew-adjoint Fredholm operators~$P$, the mod~2 dimension $\dim \ker
P\pmod2$ is a deformation invariant: skew-adjointness implies that nonzero
eigenvalues move in pairs as $P$~deforms, so in particular $\dim\ker P$ can
jump only by even numbers.\footnote{The spectrum of a Fredholm operator need
not be discrete; still, the assertion about the kernel is valid.}
For the Dirac operators~$D,D'$ on the circle, this \emph{mod~2 index}
distinguishes the two spin structures, which we have rendered here in a
concrete form. 
 
The symbol~$p(x,\xi )$ of a \emph{real} elliptic operator~$P$ on~$\RR^n$, as
a function of $\xi \in {\RR^n}^*$, is essentially the Fourier transform of a
real-valued function; see~\eqref{eq:46}.  Thus $p$~is a complex-valued
function which satisfies the reality condition
  \begin{equation}\label{eq:55}
     p(x,-\xi )=\overline{p(x,\xi )}.
  \end{equation}
The global symbol~$\sigma (P)$ on a manifold~$X$ is a map of complex vector
bundles~\eqref{eq:29}, and the reality condition~\eqref{eq:55} globalizes to
a complex conjugate isomorphism of the vector bundles which covers the
$(-1)$-involution on~$T^*X$ and commutes with~$\sigma (P)$.  This idea led
Atiyah to develop Real $K$-theory~\cite{A3}; see~\cite[\S8]{Seg4}.  Building
on this, Atiyah-Singer~\cite{AS5} proved an index theorem for real elliptic
operators.  This index theorem covers mod~2 indices, such as the invariant
on~$\cir$ mentioned above, and other cases as well.  For example, the natural
integer invariant of the real Dirac operator in dimension four is one-half
the integer index, which leads to an analytic proof of a theorem of
Rokhlin~\cite{Roh}: the signature of a closed spin 4-manifold is divisible
by~16.  These real index theorems are illuminated by Clifford algebras, as
explained in~\S\ref{subsec:5.6}.

  \begin{remark}[]\label{thm:10}
 By design, the $K$-theory proof~(\S\ref{subsec:4.4}) is suited to prove
generalizations of Theorem~\ref{thm:8}, such as the real index theorem.  For
the mod~2 indices I am unaware of other proofs, whereas for integer indices
there are other approaches, such as the heat equation methods we survey
in~\S\ref{sec:6}.
  \end{remark}

  \subsection{Index theorems for Clifford linear operators}\label{subsec:5.6}

Let $X$~be a Riemannian spin manifold of dimension~$n$, and $\Spin(X)\to X$
the principal $\Spin_n$-bundle which encodes the Riemannian spin structure;
see~\S\ref{subsec:4.1}.  Recall the Clifford algebras~$\Clpm$ with
generators~$\gamma ^1,\dots ,\gamma ^n$ and relations 
  \begin{equation}\label{eq:79}
     \gamma ^i\gamma ^j + \gamma ^j\gamma ^i = \pm 2\delta ^{ij} = \begin{cases}
     \pm 2,&i=j;\\\phantom{\pm}0,&i\neq j,\end{cases}\qquad 1\le i,j\le n. 
  \end{equation}
Use the ``mixing construction'' to form the vector bundle 
  \begin{equation}\label{eq:80}
     S_X=\Spin(X)\times \mstrut _{\Spin_n}\Cnp\longrightarrow X, 
  \end{equation}
where $\Spin_n\subset \Cnp$ acts on~$\Cnp$ by left multiplication.  The
fibers of $S_X\to X$ are $\zt$-graded \emph{right} $\Cnp$-modules, since
right and left multiplication commute.  Then since $\Cnm$~is the
($\zt$-graded) opposite algebra to~$\Cnp$, the fibers of $S_X\to X$ are
$\zt$-graded \emph{left} $\Cnm$-modules.  Furthermore, $S_X\to X$ inherits a
covariant derivative from the Levi-Civita connection.  The Dirac
operator~\eqref{eq:41} acts on sections of $S_X\to X$, where we use the
embedding $\RR^n\hookrightarrow \Cnp$ as the span of $\gamma ^1,\dots ,\gamma
^n$ to define Clifford multiplication.  Then $D_X$~is a skew-adjoint Fredholm
operator (on suitable Sobolev spaces), it is odd with respect to the grading,
and it graded commutes with the left $\Cnm$-action.

  \begin{remark}[]\label{thm:14}
\ 

 \begin{enumerate}
 \item This \emph{Clifford linear Dirac operator} appears in~\cite[\S
II.7]{LM},  attributed to Atiyah-Singer.

 \item As for any irreducible representation, the spin
representations~$\SS^{0},\SS^1$ in~\S\ref{subsec:4.1} are only determined up to
tensoring with a line, as already remarked there.  By contrast, the
construction of $S_X\to X$ and~$D_X$ are canonical for Riemannian spin
manifolds.  We have traded the irreducible representations~$\SS^0,\SS^1$ for
the canonical representation~$\Cnp$ equipped with a canonical commuting
algebra of operators.

 \item One can complexify~$\Cnp$ and so $S_X\to X$ to construct a canonical
complex Dirac operator, but the real Dirac operator contains more refined
information.

 \item If $V\to X$ is a real vector bundle with covariant derivative, we can
form a twisted Dirac operator~$D_X(V)$ on sections of $S_X\otimes V\to X$; it
too commutes with the left $\Cnm$-action.

 \end{enumerate}
  \end{remark}

The kernel and cokernel of~$D_X(V)$ are finite dimensional $\Cnm$-modules, so
by~\cite{ABS} their formal difference---the index---represents an element of
the real $K$-theory group~$KO^{-n}(\pt)$.  This is the \emph{analytic
pushforward} of~$V$.  The spin structure orients the map $f\:X\to \pt$ for
real $K$-theory---see~\eqref{eq:19} for the analogous situation in complex
$K$-theory---and the \emph{topological pushforward} is $f_![V]\in
KO^{-n}(\pt)$.  The index theorem for Dirac operators asserts the equality of
these pushforwards.  The $KO$-groups of a point are computed from Bott
periodicity.  Let $O_\infty $~be the stable orthogonal group.  Then
for~$n\ge1$,
  \begin{equation}\label{eq:81}
     KO^{-n}(\pt)\cong \pi _{n-1}O_\infty \cong \begin{cases} \ZZ,&n\equiv
     0\pmod8;\\ \zt,&n\equiv 1\pmod8;\\ \zt,&n\equiv 2\pmod8;\\ 0,&n\equiv
     3\pmod8;\\ \ZZ,&n\equiv 4\pmod8;\\ 0,&n\equiv 5\pmod8;\\ 0,&n\equiv
     6\pmod8;\\ 0,&n\equiv7\pmod8.\end{cases} 
  \end{equation}
($KO^0(\pt)\cong \ZZ$ is also correct.)  The mod~2 index for~$n=1$ discussed
in~\S\ref{subsec:5.3} is a special case. 
 
In~\cite{AS6}, written in 1969, Atiyah-Singer study spaces~$\Fn(\sH_n)$,
$n\in \ZZ$, of odd skew-adjoint Fredholm operators acting on real and complex
$\zt$-graded Hilbert spaces~$\sH_n$ equipped with a $\Cn$~action; the
operators graded commute with~$\Cn$.  The Clifford linear Dirac operator is
an example.  The main result is that the spaces~$\Fn(\sH_n)$ form a spectrum
in the sense of homotopy theory, and Atiyah-Singer identify it with the
$KO$-~and $K$-theory spectra in the real and complex cases, respectively.
This generalizes the relationship between $K$-theory and Fredholm operators
in earlier work of J\"anich~\cite{J}.  It also gives an alternative
construction of the analytic index of the Clifford linear Dirac operator.
Finally, the periodicity of Clifford algebras~\cite{ABS} leads to another
proof of Bott periodicity; see also~\cite{Kar,Wo}.

  \subsection{Families index theorem}\label{subsec:5.4}

Recall that Grothendieck, in his Riemann Roch Theorem~\ref{thm:4}, considers
proper morphisms $f\:X\to S$, not just a single variety $X\to\pt$.
Similarly, Atiyah-Singer~\cite{AS4} generalize Theorem~\ref{thm:8} to proper
fiber bundles $f\:X\to S$ equipped with a family~$P$ of elliptic
(pseudo)differential operators along the fibers.  On suitable function
spaces, $P$~is a continuous family of Fredholm operators parametrized by~$S$,
and the analytic index is the homotopy class of this family.  Since Fredholm
operators form a classifying space for $K$-theory~\cite{J,A4,AS6}, the
analytic index lies in~$K^{\bullet }(S)$.  The topological index
construction, executed in the fiber bundle~$f$, also leads to an element
of~$K^{\bullet }(S)$.  The \emph{index theorem for families} is the equality
of the analytic and topological indices.
 
There is also a families version of the Clifford linear story
of~\S\ref{subsec:5.6}, which we illustrate.

  \begin{example}[]\label{thm:24}
 Let $f\:X\to S$ be a proper fiber bundle with fibers of odd dimension~$n$.
Suppose a relative Riemannian spin structure\footnote{A \emph{relative spin
structure} is a spin structure on the vertical tangent bundle $T(X/S)\to X$.
A \emph{relative Riemannian structure} is an inner product on this vector
bundle together with a horizontal distribution on the fiber bundle~ $f$.} is also
given.  In \emph{complex} $K$-theory there is a pushforward
  \begin{equation}\label{eq:82}
     f_!\:K^0(X)\longrightarrow K^{-n}(S).
  \end{equation}
If $V\to X$ is a complex vector bundle with covariant derivative, we form the
family~$D_{X/S}(V)$ of complex Clifford linear Dirac operators parametrized
by~$S$.  The analytic index is the homotopy class of the map $S\to\Fred_{-n}$
given by the Dirac operators, where $\Fred_{-n}$ is the space of Fredholm
operators introduced at the end of~\S\ref{subsec:5.6}, a classifying space
for~$K^{-n}$.  The index theorem asserts that the analytic index equals
$f_![V]$.  The ``lowest piece'' of the index is captured by composing with a
natural map 
  \begin{equation}\label{eq:83}
     K^{-n}(S)\longrightarrow H^1(S;\ZZ), 
  \end{equation}
which may be considered a determinant map.  Recalling that a class
in~$H^1(S;\ZZ)$ is determined by its periods on maps $\varphi \:\cir\to S$,
in other words by its values on $\varphi _*[\cir]\in H_1(S)$, since
$H^1(S;\ZZ)$ is torsionfree, we compute the image of $\ind D_{X/S}(V)$
under~\eqref{eq:83} by base changing via~$\varphi $ to a fiber bundle
$X_\varphi \to\cir$ equipped with a vector bundle $V_\varphi \to X_\varphi $.
Then $X_\varphi $~is an even dimensional closed Riemannian spin manifold, and
the integer period we seek is the numerical index $\ind D\mstrut _{X_\varphi
}(V_\varphi )\in \ZZ$.
  \end{example}

  \begin{remark}[]\label{thm:15}
 \ 
 \begin{enumerate}

 \item The family of Clifford linear Dirac operators $D_{X/S}(V)$ gives rise
to a family of self-adjoint Dirac operators on an ungraded Hilbert space, a
formulation which appears more frequently; see \cite[\S3]{APS3}, for example.

 \item A cohomology class in~$H^1(S;\ZZ)$ is a homotopy class of maps $S\to
\RZ$.  A \emph{geometric} invariant of the family~$D_{X/S}(V)$---the
Atiyah-Patodi-Singer $\eta $-invariant---promotes this homotopy class to a
specific map, as we will see in~\S\ref{subsec:7.3}.

 \item Although we stated these topological constructions for Dirac
operators, they generalize in various ways to families of elliptic
pseudodifferential operators.

 \end{enumerate}
  \end{remark}

A situation in which one encounters a family of \emph{linear} elliptic
equations is linearization of solutions to a \emph{nonlinear} elliptic
equation.  The abstract setup is a nonlinear Fredholm map $F\:\sB\to\sC$
between infinite dimensional Hilbert manifolds.  Sometimes $F$~is equivariant
for the action of an infinite dimensional Lie group~$\sG$ and is only
Fredholm modulo~$\sG$.  One is then interested in the moduli space, or
stack,~$\sM=F\inv (c)/\sG$ for some~$c\in \sC$.  If $F$~is a nonlinear
elliptic operator on a closed manifold, mapping between functions spaces,
then its linearizations~$dF_b$ at~$b\in F\inv (c)$ fit together to a family
of linear elliptic operators parametrized by~$\sB$.  The index of this family
computes the (virtual) tangent bundle to~$\sM$ and yields useful information
about~$\sM$.  This general plan is used by Atiyah-Hitchin-Singer~\cite{AHS}
to investigate the instanton equations on a 4-manifold (see Donaldson's
paper~\cite[\S2]{Do} in this volume) and it has been used since in many other
problems in geometric analysis.

  \begin{remark}[]\label{thm:11}
 The image of the families index under the Chern character $\ch\:K^{\bullet
}(S)\to H^{\bullet }(S;\QQ)$ is a cruder invariant than the $K$-theory index,
but often it contains information of interest.  It can be computed by a
topological analogue of~\eqref{eq:12}, and it is accessible via heat equation
methods~\cite{Bi2}, whereas the more powerful $K$-theory index is not, as far
as we know. 
  \end{remark}

The Atiyah-Hitchin-Singer work is the first of many applications of the
families index theorem to quantum field theory and string theory, and to
mathematical problems arising from that physics.  We take up an additional
example in~\S\ref{sec:8}.

  \subsection{Coverings and von Neumann algebras}\label{subsec:5.5}

As should be clear by now, Atiyah used the index theorem as a launching pad
for mathematical adventures in many directions.  One which proved
particularly fruitful involves von Neumann algebras~\cite{A5}.
Atiyah~\cite{A1} comments:
 
  \begin{quote}
 In particular I learnt from Singer, who had a strong background in
functional analysis, about von Neumann algebras of type~II with their
peculiar real-valued dimensions.  We realized that $K$-theory and index
theory could be generalized in this direction, but it was not clear at first
if such a generalization would really be of any interest.  However in one
particularly simple case, that of a manifold with an infinite fundamental
group, it became clear that the ideas of von Neumann algebras were quite
natural and led to concrete non-trivial results.  This was the content of my
talk~\cite{A5} at the meeting in honour of Henri Cartan.  Since I was not an
expert on von Neumann algebras I attempted in this presentation to give a
simple, elementary and essentially self-contained treatment of the results.
Later on in the hands of Alain Connes, the world expert on the subject, these
simple ideas were enormously extended and developed into a whole theory of
linear analysis for foliations.
  \end{quote}

\noindent 
 And in the hands of Connes, Kasparov and many others into index theory and
$K$-theory for $C^*$-algebras.  An influential conference talk indeed!

The situation in~\cite{A5} is an unramified Galois covering $\pi \:\tX\to X$
with Galois group~$\Gamma $ acting freely on~$\tX$, and a $\Gamma $-invariant
elliptic operator~$\tD$ on~$\tX$.  There is an induced elliptic operator~$D$
on~$X$, and we assume $X$~is compact.  For example, $X$~could be a closed
Riemann surface and $\tX$~its universal cover which, if the genus of~$X$ is
~$\ge2$, is isomorphic to the unit disk~$\Omega $.  If $\tD$~is the
$\dbar$-operator, then $\ker\tD$~is the \emph{infinite dimensional} space of
holomorphic functions on~$\Omega $.  In general, if $\Gamma $~is infinite
then $\ker\tD$ and $\coker\tD$~are infinite dimensional, whereas $\ker D$~and
$\coker D$~are finite dimensional, the latter since $X$~is assumed compact.
So $\ind D$~is well-defined.  Atiyah introduces a $\Gamma $~invariant
measure on~$\tX$ and the von Neumann algebra~$\ssA$ of bounded linear
operators on $L^2(\tX,\tE^0)$ for $\tE^0\to\tX$ the vector bundle on whose
sections $\tD$~is defined.  Orthogonal projection onto~$\ker\tD$ lies in the
von Neumann algebra, and its von Neumann trace, a real number, is defined to
be the $\Gamma $~dimension of~$\ker\tD$.  Repeating for~$\coker\tD$, Atiyah
defines a real-valued index $\ind_\Gamma \tD$. 

  \begin{theorem}[Atiyah~\cite{A5}]\label{thm:12}
 $\ind_\Gamma \tD = \ind D$. 
  \end{theorem}

\noindent
 Atiyah's account of this theorem, as stated earlier, was a catalyst for
index theory on noncompact spaces, singular spaces, and beyond.

   \section{Heat equation proof}\label{sec:6}

Beginning in the late~1960's the expanding circle of ideas emanating from the
basic Atiyah-Singer index theorem took a more analytic turn.  The focus
shifted beyond the kernel of elliptic operators to include higher
eigenvalues.  This led first to a \emph{local} version of the index theorem
and then to local \emph{geometric} invariants (as opposed to \emph{global
topological} invariants).  We treat the former in this section and the latter
in the next.
 
We begin in~\S\ref{subsec:6.1} with two basic constructs to collate higher
eigenvalues into a single function: the $\zeta $-function and the trace of
the heat operator.  (They are analogues of basic objects in analytic number
theory.)  The local index theorem, proved first in special cases by Vijay
Patodi and then in general by Peter Gilkey, is the subject
of~\S\ref{subsec:6.2}. This work was completed in the early 1970's.  In the
1980's several new proofs of the local index theorem led to a deeper
understanding of the origins of the $\Ahat$-genus in the index formula for
Dirac operators.  In~\S\ref{subsec:6.3} we briefly summarize these
contributions by Ezra Getzler, Edward Witten, Jean-Michel Bismut, Nicole
Berline, Michele Vergne, and of course Michael Atiyah.

  \subsection{Heat operators, zeta functions, and the index}\label{subsec:6.1}

Let $\Delta $~be a nonnegative self-adjoint operator on a Hilbert
space~$\sH$.  We seek to define the \emph{heat operator}
  \begin{equation}\label{eq:56}
     H_t=e^{-t\Delta },\qquad t\in \RR^{>0}, 
  \end{equation}
and the \emph{$\zeta $-function}
  \begin{equation}\label{eq:57}
     \zD(s)=\Tr\Delta ^{-s},\qquad s\in \CC. 
  \end{equation}
Both are well-defined if $\sH$~is finite dimensional, and they are related by
the Mellin transform 
  \begin{equation}\label{eq:58}
     \Tr\Delta ^{-s} = \frac{1}{\Gamma
     (s)}\int_{0}^{\infty}\frac{dt}{t}\,t^s\,\Tr\bigl(e^{-t\Delta } \bigr). 
  \end{equation}
 
If $\Delta $~is a nonnegative self-adjoint second-order elliptic operator on
a closed manifold~$X$, then the heat operator~$H_t$ in~\eqref{eq:56} is
exists by basic elliptic theory.  It is a \emph{smoothing operator}: $H_t$~
maps distributions to smooth functions.  For example, if $\delta _y$~is the
Dirac $\delta $-distribution at~$y\in X$, then\footnote{Our notation assumes
$\Delta $~is an operator acting on functions.  A small modification
incorporates vector bundles.}
  \begin{equation}\label{eq:59}
     h_t(x,y) = \left( e^{-t\Delta }\delta _y \right)\!(x),\qquad x\in X, 
  \end{equation}
is a smooth function of~$t,x,y$ called the \emph{heat kernel}.  If $\Delta
$~is the scalar Laplace operator attached to a Riemannian metric on~$X$---the
Laplace-Beltrami operator---then intuitively $h_t(x,y)$~is the amount of heat
at~$x$ after time~$t$ given an initial distribution~$\delta _y$ of heat.
Heat flows with infinite propagation speed and instantly diffuses:
$h_t(x,y)>0$ for all~$t,x,y$.  The properties of physical heat flow inform
intuition about the large and small time behavior of the heat operator of a
general nonnegative self-adjoint second-order differential operator~$\Delta $
acting on sections of a vector bundle $E\to X$.  As $t\to\infty $ the heat
operator~$e^{-t\Delta }$ converges (in the uniform topology) to projection
onto~$\ker \Delta $.  As $t\to0$ the heat operator converges (in the strong
operator topology) to the identity operator.  A more precise version of the
small~$t$ behavior is the subject of an influential 1948 paper of
Minakshisundarum and Pleijel~\cite{MP} in case $\Delta $~is the
Laplace-Beltrami operator; Seeley~\cite{Se4} extends their results to more
general elliptic pseudodifferential operators.  As $t\to0$ the heat
kernel~$h_t(x,y)$ converges exponentially to zero if~$x\neq y$, and on the
diagonal there is an asymptotic expansion
  \begin{equation}\label{eq:60}
     h_t(x,x)\sim t^{-n/2}\sum\limits_{k=0}^{\infty}A_k(x)t^i\qquad
     \textnormal{as $t\to0$}, 
  \end{equation}
where $A_k$~are smooth functions on~$X$.  For~$x\in X$, the value of
$A_k(x)\in \End E_x$ depends only on a finite jet of the total symbol of the
differential operator~$\Delta $ at~$x$.  If $\Delta $~is canonically
associated to a Riemannian metric, then $A_k(x)$~depends on a finite jet of
the metric at~$x$; the order of the jet grows with~$k$.
 
For elliptic operators~$\Delta $ on compact manifolds of the type discussed
in the previous paragraph, the $\zeta $-function~\eqref{eq:57} exists and
$\zD$~is a holomorphic function of~$s$ for $\Re (s)>\!>0$.  The asymptotic
expansion~\eqref{eq:60} of the heat kernel~$h_t$ is equivalent, via the
Mellin transform~\eqref{eq:58}, to a meromorphic extension of~$\zD$ to the
entire complex $s$-line, which in fact is what is proved in~\cite{MP,Se4}.
 
Now suppose $P\:C^{\infty}(X,E^0)\to C^{\infty}(X,E^1)$ is a first-order elliptic
operator on a closed $n$-manifold~$X$ equipped with complex vector bundles
$E^0,E^1\to X$.  Assume metrics everywhere so that the formal adjoint
$P^*\:C^{\infty}(X,E^1)\to C^{\infty}(X,E^0)$ is defined; it too is a first-order
elliptic operator.  Then 
  \begin{equation}\label{eq:61}
     \ind P = \dim\ker P^*P - \dim\ker PP^*, 
  \end{equation}
and each of $P^*P,\,PP^*$ is a nonnegative self-adjoint second-order elliptic
differential operator.  For~$\lambda \ge0$, let $\sE^0_\lambda \subset
C^{\infty}(X,E^0)$ be the $\lambda $-eigenspace of~$P^*P$ and $\sE^1_\lambda
\subset C^{\infty}(X,E^1)$ the $\lambda $-eigenspace of~$PP^*$.  Then
for~$\lambda >0$, 
  \begin{equation}\label{eq:62}
     P\res{\sE^0_\lambda }\:\sE^0_\lambda \longrightarrow \sE^1_\lambda 
  \end{equation}
is an isomorphism.  Therefore, for any function $\chi
\:\RR^{\ge0}\to\RR^{\ge0}$ such that $\chi (0)=1$ and $\chi (\lambda
)\searrow0$ sufficiently rapidly as $\lambda \to\infty $, 
  \begin{equation}\label{eq:63}
     \ind P = \sum\limits_{\lambda \in \spec P^*P}\chi (\lambda ) -
     \sum\limits_{\lambda \in \spec PP^*}\chi (\lambda ). 
  \end{equation}
For $\chi (\lambda )=\lambda ^{-s}$ with $\Re(s)>\!>0$, we obtain a formula
for the index which appears in Atiyah-Bott~\cite[\S8]{AB1}:
  \begin{equation}\label{eq:64}
     \ind P = \Tr \zeta \mstrut _{P^*P}(s) - \Tr \zeta \mstrut _{PP^*}(s). 
  \end{equation}
In fact \eqref{eq:64}~holds for all~$s\in \CC$, due to the meromorphic
continuation of zeta functions.  Atiyah-Bott note that~$s=0$ is a
particularly good argument in view of explicit integral formulas~\cite{Se4}
in terms of the symbol of~$P$.  (Another motivation for setting~$s=0$: for an
operator $P\:\sH^0\to\sH^1$ between finite dimensional Hilbert spaces, the
value at~$s=0$ is $\dim \sH^0 - \dim \sH^1$.)  But while the explicit
formulas are local, they involve high derivatives of the symbol, whereas the
characteristic class formula~\eqref{eq:44} for the index only involves a few
derivatives when written in terms of Chern-Weil polynomials of the
curvature.  It is this mismatch which remained a mystery for several years. 
 
In place of $\zeta $-functions, the trace of the heat kernel is commonly used
in~\eqref{eq:63}.  This corresponds to $\chi (\lambda
)= e^{-t\lambda }$, $t>0$.  Then for all~$t>0$, we have
  \begin{equation}\label{eq:65}
     \ind P = \Tr e^{-tP^*P} - \Tr e^{-tPP^*}. 
  \end{equation}
In fact, one can prove the right hand side is constant in~$t$ by
differentiation, and evaluation as~$t\to\infty $ reproduces~\eqref{eq:65}.
On the other hand, let~$t\to0$ and use the asymptotic expansion~\eqref{eq:60}
to obtain 
  \begin{equation}\label{eq:66}
     \ind P = \int_{X}\tr\left[ A^0_{n/2}(x) - A^1_{n/2}(x) \right] \,|dx|, 
  \end{equation}
where $A^0_k(x),A^1_k(x)$~are the heat coefficients for~$P^*P,PP^*$ acting on
the vector spaces~$E^0_x,E^1_x$, respectively.  As is true for the $\zeta
$-function, the formulas for~$A^i_{n/2}(x)$ involve many derivatives of the
symbol of~$P$ at~$x$, so seem inaccessible as a means of proving the index
formula.

  \subsection{The local index theorem}\label{subsec:6.2}

For the Laplace-Beltrami operator~$\Delta $ on an $n$-dimensional Riemannian
manifold~$X$, the first coefficient~$A_0$ in~\eqref{eq:60} is the constant
function $1/(4\pi )^{n/2}$.  This reflects the solution to the classical heat
equation in Euclidean space, and it implies Weyl's law for the asymptotic
growth of the eigenvalues of~$\Delta $, which only depends on~$n$ and
$\Vol(X)\sim\int_{X}A_0(x)\,|dx|$.  Weyl's law, which does not depend on the
heat kernel expansion, was one motivation for Mark Kac~\cite{K} to ask
in~1966: To what extent do the eigenvalues of~$\Delta $ determine the
Riemannian manifold~$X$?  Kac focused on domains in the Euclidean
plane~$\EE^2$, though the more general question is implicit.
McKean-Singer~\cite{MS} immediately took this up, and they determined the
next few coefficients~$A_1,A_2$ in the heat kernel expansion, thereby proving
a conjecture of Kac-Pleijel.  (There are contemporaneous independent results
by de Bruijn, Arnold, and Berger.)  In particular, $A_1$ is a multiple of the
scalar curvature.  For~$n=2$, McKean-Singer observe a cancellation which
holds at each point of~$X$, and they conjecture a similar result in all
dimensions.  Namely, let $\Delta \q$~denote the Laplace operator on
differential $q$-forms, $h_t\q$~the associated heat kernel, and $A_k\q$ the
heat coefficients.  The McKean-Singer conjecture is that for all~$x\in X$ the
limit
  \begin{equation}\label{eq:67}
     \lim\limits_{t\to0}\;\sum\limits_{q=0}^n\,(-1)^q\tr h_t\q(x,x) 
  \end{equation}
exists and, furthermore, for $n$~even it equals the Gauss-Bonnet-Chern
integrand which integrates to the Euler number of~$X$.  From~\eqref{eq:60}
the existence of the limit is equivalent to the cancellation  
  \begin{equation}\label{eq:68}
     \sum\limits_{q=0}^n(-1)^q\tr A_k\q(x)=0,\qquad k<\frac n2,\quad x\in
     X. 
  \end{equation}
The alternating sum for~$k=n/2$ equals the limit~\eqref{eq:67}.  McKean-Singer
prove that the limit exists and vanishes for $n$~odd, and they compute the
limit for~$n=2$.   
 
In~1970 Patodi~\cite{P1} proved the McKean-Singer conjecture via a virtuoso
direct computation.  He immediately~\cite{P2} applied his methods to prove
the Riemann-Roch Theorem~\ref{thm:1} for K\"ahler manifolds.  A few years
later, Gilkey used different methods in his PhD thesis~\cite{Gi1}---scaling
plays a crucial role---and proved the corresponding theorem for
twisted signature operators.  By standard topological arguments this implies
the Atiyah-Singer Index Theorem~\ref{thm:5}; see~\S\ref{subsec:4.2}.
Subsequently, Atiyah-Bott-Patodi~\cite{ABP} gave a proof of Gilkey's Theorem
and of the resulting proof of the index theorem.  We recount Gilkey's main
result.
 
Gilkey investigates differential forms built canonically from a Riemannian
metric.  Using modern terminology to economize, let $\Mn$~be the category of
smooth $n$-manifolds and local diffeomorphisms.  Consider the functors
(sheaves)
  \begin{equation}\label{eq:69}
     \begin{aligned} \Met\:\Mn\op&\longrightarrow \Set \\ \Omega
      ^q\:\Mn\op&\longrightarrow \Set\end{aligned} 
  \end{equation}
where if $M$~is a smooth $n$-manifold, then $\Met(M)$~is the set of
Riemannian metrics on~$M$ and $\Omega ^q(M)$~is the set of differential
$q$-forms.  We seek natural transformations 
  \begin{equation}\label{eq:70}
     \omega \:\Met\longrightarrow \Omega ^q. 
  \end{equation}
Roughly speaking, these are assignments of differential forms to Riemannian
metrics covariant under coordinate changes.  More poetically, they are
differential forms on~$\Met$.  Even for~$q=0$ the classification problem is
intractable: any smooth function of the scalar curvature gives a natural
function of the Riemannian metric.  Now introduce scaling.  We say $\omega
$~is \emph{homogeneous of weight~$k$} if
  \begin{equation}\label{eq:71}
     \omega (\lambda ^2g)=\lambda ^k\omega (g)\qquad \textnormal{for all
     $\lambda \in \RR^{>0}$}. 
  \end{equation}
We say $\omega $~is \emph{regular} if in any local coordinate system
$x^1,\dots ,x^n$ it takes the form 
  \begin{equation}\label{eq:72}
     \omega (g)(x) = \sum\limits_{I}\sum\limits_{\alpha
     }^{\textnormal{finite}}\sum\limits_{i,j=1}^n \omega ^{i,j}_{I,\alpha
     }(x)\, \frac{\partial ^{|\alpha |}g_{ij}}{\partial x^{\alpha
     _1}\cdots\partial x^{\alpha _n}} \,dx^{i_1}\wedge \cdots\wedge
     dx^{i_q}, 
  \end{equation}
where $I=(i_1,\dots ,i_q)$ with $1\le i_1<\cdots< i_q\le n$, and $\alpha
=(\alpha _1,\dots ,\alpha _n)$ with $\alpha _k\in \ZZ^{\ge0}$.  The $\omega
^{ij}_{I,\alpha }$~are smooth functions.  Crucially, only a finite set
of~$\alpha $ appears.

  \begin{theorem}[Gilkey~\cite{Gi1}]\label{thm:13}
 A natural differential form~\eqref{eq:70} which is regular and homogeneous
of nonnegative weight is a polynomial in the Chern-Weil forms of the
Pontrjagin classes. 
  \end{theorem}

\noindent
 The nonzero forms have weight zero. 
 
The proof of Gilkey's Theorem~\ref{thm:13} in~\cite{ABP} uses Weyl's
theorem~\cite{W} on invariants of the orthogonal group.  Atiyah-Bott-Patodi
apply Theorem~\ref{thm:13} to the signature operator~$P$ on a Riemannian
manifold~$X$.  Resuming the notation of~\S\ref{subsec:6.1}, the vanishing of
positive weight forms implies
  \begin{equation}\label{eq:73}
     \tr\bigl[A^0_k(x) - A^1_k(x) \bigr]=0,\qquad k<\frac n2,\quad x\in X. 
  \end{equation}
This cancellation result implies, as in~\eqref{eq:68}, the existence of 
  \begin{equation}\label{eq:74}
     \lim\limits_{t\to0}\;\bigl[\tr h^0_t(x,x) - \tr h^1_t(x,x) \bigr], 
  \end{equation}
and Theorem~\ref{thm:13} tells that the limit is a polynomial in Pontrjagin
forms.  The precise polynomial---the $L$-genus~\eqref{eq:9}---is determined
as in Hirzebruch's original proof by computing enough examples.  The road
from here to the \emph{global} index theorem follows established lines.  It
is the \emph{local} index theorem---the existence and identification of the
limit~\eqref{eq:74}---which leads to future developments.

  \subsection{Postscript: Whence the $\Ahat$-genus?}\label{subsec:6.3}

Different conceptual understandings of the cancellation~\eqref{eq:73} and of
the limiting value~\eqref{eq:74} were achieved in the first half of the
1980's.  The setting is (generalized) Dirac operators, where basic properties
of Clifford algebras yield the cancellation.  The limit is the Chern-Weil
$\Ahat$-form, a polynomial in the Pontrjagin forms, whose appearance is
derived from various sources.  In these works the $\Ahat$-genus appears by
direct argument.  We give a brief resum\'e.
 
One route to the $\Ahat$-genus passes through Mehler's formula for the heat
kernel of the harmonic oscillator~\cite[p.~19]{GJ}, which Getzler~\cite{Ge2}
employs in his proof of the local index theorem.  He uses a homothety which
not only scales time and space, but also scales the Clifford algebra
variables in the Dirac operator~\eqref{eq:41}.  His technique was in part
inspired by contemporary physics proofs of the index
theorem~\cite{Wi1,Ag,FW,Ge1} using supersymmetric quantum mechanics.

At a conference in honor of Laurent Schwartz, Atiyah~\cite{A6} exposed
Witten's idea to derive the index theorem by applying the Duistermaat-Heckman
exactness of stationary phase theorem~\cite{DH} to the free loop space of a
compact Riemannian manifold.  (During that period Atiyah-Bott~\cite{AB5}
placed the Duistermaat-Heckman result in the context of localization in
equivariant cohomology.)  In this proof the $\Ahat$-genus enters by
regularizing a certain infinite product, as it does in the supersymmetric
quantum mechanics proof~\cite[\S1.2.4]{Wi2}.
 
Inspired by Atiyah's account, Bismut~\cite{Bi1} executed a proof of the index
theorem using Wiener measure on loop space and Malliavin calculus.  In this
way he deals with integrals over loop space rigorously.  The heat kernel is
represented in terms of Wiener measure with the aid of Lichnerowicz's
formula, which expresses the Dirac Laplacian in terms of the covariant
Laplacian.  The localization to point loops as~$t\to0$ is natural in this
probabilistic approach.  The variable~$t$ represents the total time during
which a Brownian path exists, and as the time tends to zero, only constant
loops have a significant probability of occurring.  The evaluation of the
integral over these point loops is accomplished using a formula of Paul
L\'evy~\cite{L}, who considers a Brownian curve in the plane conditioned to
close after time~$2\pi $.  Then the characteristic function of the area~$S$
enclosed by the random curve (expectation value of~$e^{izS}$, $z\in \CC$)
is~$\pi z/\sinh \pi z$.  This same calculation appears in Bismut's work, only
there the curvature of~$X$ replaces~$z$, and once again the $\Ahat$-genus is
obtained.

The $\Ahat$-genus arises quite differently in a proof of the index theorem
due to Berline and Vergne~\cite{BV}.  Let $G$~be a Lie group with Lie
algebra~$\frak{g}$.  Then a standard formula in the theory of Lie groups
asserts that the differential of the exponential map~$\exp\:\frak{g}\to G$
at~$a\in \frak{g}$ is
  \begin{equation}\label{eq:76}
     d\exp_a = \frac{1-e^{-\ad a}}{\ad a} = J (\ad a),
  \end{equation}
where the power series which defines~$J$ is the multiplicative inverse of the
power series which defines the Todd genus~\eqref{eq:5}.  It was a mystery
whether the occurrence of the Todd genus in~\eqref{eq:76} is related to the
index theorem.  Berline and Vergne noticed that if $X$~is a Riemannian
manifold, and $\O(X)$~the principal bundle of orthonormal frames, then the
differential of the Riemannian exponential map on~$\O(X)$ is given by a
similar formula.  Precisely, there is a natural isomorphism $T_p\O(X)\cong
\Rn\oplus\frak{o}(n) $ via the Levi-Civita connection, and the differential
of the exponential map $\exp\:\Rn\oplus\frak{o}(n)\to \O(X)$ at~$p\in \O(X)$,
evaluated on~$a\in \frak{o}(n)$, is
  \begin{equation}\label{eq:75}
     \begin{aligned} &d\exp_a\res{\Rn} = \exp(-a)\,J \bigl( \langle
     \Omega _p/2,a\rangle\bigr).\\ &d\exp_a\res{\frak{o}(n)}=J
     (\ad a); 
     \end{aligned} 
  \end{equation}
In this formula the Riemann curvature~$\Omega $, which takes values
in~$\frak{o}(n)$, is contracted with~$a$ using the Killing form.  The result
is a 2-form, which can be identified as an element of~$\frak{o}(n)$.  To
prove the index theorem Berline and Vergne work on the frame bundle~$\O(X)$,
not on the base~$X$.  To compensate for the introduction of extra degrees of
freedom in the fiber direction, they must study the behavior of the heat
kernel along the fiber.  It is at this stage, in the small time limit, where
\eqref{eq:75}~appears.  Ultimately, that is how the $\Ahat$-genus enters
their proof.

   \section{Geometric invariants of Dirac operators}\label{sec:7}

Up to this point index theory produced \emph{global topological} invariants
of elliptic operators; their natural home is topological $K$-theory.
Beginning with work of Atiyah-Patodi-Singer announced in~1973, index theory
took a turn towards \emph{local differential geometric}
invariants.\footnote{The natural home for the geometric invariants is
\emph{differential} $K$-theory, but that is a more recent development beyond
the scope of this article.}  Furthermore, the focus shifted from general
elliptic pseudodifferential operators to Dirac operators.  Heat equation
methods provide the fundamental tools to construct invariants. 
 
The three papers~\cite{APS1,APS2,APS3} of Atiyah-Patodi-Singer introduce the
$\eta $-invariant of a Dirac operator.  Its definition is parallel to that of
an $L$-function in analytic number theory.  Their first main theorem, which
we recount in~\S\ref{subsec:7.3}, is an index theorem for a Dirac operator on
a compact Riemannian spin manifold with boundary.  A key ingredient in the
story are new \emph{global} elliptic boundary
conditions~(\S\ref{subsec:7.2}); local elliptic boundary conditions are
obstructed in most cases, as Atiyah-Bott had discovered a decade earlier.
The $\eta $-invariant solves a problem which served as one motivation for
their work, namely the computation of the signature defect, and this is our
point of departure in~\S\ref{subsec:7.1}.  The Atiyah-Patodi-Singer papers
contain many more important theorems, such as the index theorem for flat
bundles, which we do not cover here.
 
Another view of the $\eta $-invariant is the subject of~\S\ref{subsec:7.4}.
Characteristic numbers of vector bundles over closed oriented manifolds---the
integers obtained by pairing products of Chern and Pontrjagin classes of a
vector bundle with the fundamental class of the base manifold---are primary
integer-valued topological invariants.  The associated $\RZ$-valued secondary
differential geometric invariants had been introduced by Chern and Simons a
few years prior.  Similarly, integer-valued $K$-theory characteristic
numbers, which by the index theorem are indices of Dirac operators, are
primary topological invariants.  The associated secondary differential
geometric quantity is the Atiyah-Patodi-Singer $\eta $-invariant.  The next
geometric invariant of a Dirac operator, or family of Dirac operators, is the
determinant.  The underlying theory was developed in the 1980's, as we
recount in~\S\ref{subsec:7.5}.  It is an important ingredient in the
application to physics we take up in~\S\ref{sec:8}.

  \subsection{The signature defect}\label{subsec:7.1}

Recall the classical Gauss-Bonnet theorem.  Let $X$~be a closed Riemannian
2-manifold and $K\:X\to\RR$ its Gauss curvature.  Then the Euler number
of~$X$ is the curvature integral 
  \begin{equation}\label{eq:84}
     \Euler(X) = \int_{X}\frac{K}{2\pi }\,d\mu _X, 
  \end{equation}
where $d\mu _X$~is the Riemannian measure.  If now $X$~is compact with
boundary, then there is a boundary contribution from the geodesic curvature
$\kappa \:\bX\to\RR$, namely 
  \begin{equation}\label{eq:85}
     \Euler(X) = \int_{X}\frac{K}{2\pi }\,d\mu _X + \int_{\bX}\frac{\kappa
     }{2\pi }\,d\mu _{\bX}. 
  \end{equation}
If a neighborhood of~$\bX$ in~$X$ is isometric to the cylinder $[0,\epsilon
)\times \bX$ with its product metric for some~$\epsilon >0$, then the
boundary term vanishes. 
 
Now let $X$~be a closed oriented Riemannian 4-manifold.  The Signature
Theorem~\ref{thm:3} implies 
  \begin{equation}\label{eq:86}
     \Sign(X) = \int_{X}\omega , 
  \end{equation}
where $\omega $~is the Chern-Weil 4-form of the rational characteristic
class~$p_1/3$.  If $X$~is compact with boundary, and even if we assume
the Riemannian metric is a product near the boundary, which we do, 
formula~\eqref{eq:86} need not hold.  Set~$Y=\bX$.  Then the \emph{signature
defect}~\cite[\S10.3]{A7} 
  \begin{equation}\label{eq:87}
     \alpha (Y)=\Sign(X) - \int_{X}\omega 
  \end{equation}
depends only on the closed oriented Riemannian 3-manifold~$Y$, as follows
easily from~\eqref{eq:86}.  Atiyah-Patodi-Singer~\cite{APS1} argue that
$\alpha $~is a smooth function of the Riemannian metric, is odd under
orientation-reversal, and is \emph{not} of the form $\int_{Y}\eta $ for some
natural 3-form~$\eta $ in the metric  (since the signature defect is not
multiplicative under finite covers).
 
A concrete instance of the signature defect studied by Hirzebruch~\cite{H4}
in the early 1970s was a prime motivation for Atiyah-Patodi-Singer.  Let
$K=\QQ(\sqrt d)$ be a real quadratic number field---$d\in \ZZ^{>1}$ is
assumed square-free---and let $\sO\subset K$ be the ring of integers.  The
two square roots of~$d$ give two embeddings $K\hookrightarrow \RR$, thus an
embedding
  \begin{equation}\label{eq:88}
     G=\PSL_2(\sO)\hooklongrightarrow \PSL_2(\RR)\times \PSL_2(\RR). 
  \end{equation}
Let $\HH$~be the upper half plane.  The quotient 
  \begin{equation}\label{eq:89}
     X^0(K) = \bigl(\HH\times \HH\bigr)\bigm/ G 
  \end{equation}
is a \emph{Hilbert modular surface}.  The group~$G$ may act with finite
stabilizers, i.e., $X^0(K)$~may be an orbifold which is not a smooth
manifold.  It is noncompact; $X^0(K)$~has a finite set of ends in bijection
with the ideal class group of~$K$.  Truncate each end to construct a compact
orbifold with boundary; each boundary component is a fiber bundle with
base~$\cir$ and fiber $\cir\times \cir$.  Also, cut out a neighborhood of
each orbifold point to obtain a smooth compact manifold~$X(K)$ with boundary;
the additional boundary components are lens spaces.  The signature
defect~\eqref{eq:87} at lens space boundaries was known~\cite[\S10.3]{A7}
from the equivariant signature theorem~\S\ref{subsec:5.1}.
Hirzebruch~\cite[p.~222]{H4} computed the signature defect at the other
boundaries via a desingularization of the cusp singularity in the cone on the
boundary.  The formula is a simple expression in terms of a continued
fraction associated to the singularity.  On the other hand, Shimizu~\cite{Sh}
introduced an $L$-function associated to the ideal class which corresponds to
the end.  Hirzebruch~\cite[p.~231]{H4} proves that the signature defect at
the cusp is the value of that $L$-function at~$s=1$, up to a numerical
factor.
 
There is a generalization of this story to totally real number fields of
arbitrary degree.  For the general case Hirzebruch conjectured that the
signature defect at a cusp singularity is again a value of the Shimizu
$L$-function.  This conjecture was proved independently by
Atiyah-Donnelly-Singer~\cite{ADS} and by M\"uller~\cite{Mu} in~1982--3.

  \subsection{Global boundary conditions}\label{subsec:7.2}

To compute the signature defect~\eqref{eq:87} in general, it is natural to
consider the signature operator~\eqref{eq:78} on a compact manifold with
boundary.  But we must impose an elliptic boundary condition, and for the
signature operator the topological obstruction to \emph{local} elliptic
boundary conditions~(\S\ref{subsec:5.2}) is nonzero.  Atiyah-Patodi-Singer
overcome this obstruction by a novel maneuver: they introduce \emph{global}
elliptic boundary conditions that exist for any generalized Dirac operator,
including the signature operator.  These global boundary conditions are now
ubiquitous in the theory and applications of Dirac operators on manifolds
with boundary.
 
As a first example, consider the $\dbar$-operator on the closure of the unit
disk $\Omega \subset \CC_z$, as in~\S\ref{subsec:5.2}.  The kernel consists
of holomorphic functions on~$\overline\Omega $; a dense subspace is the space
of polynomials, the linear span of $\{z^n:n\in \ZZ^{\ge0}\}$.  This is an
infinite dimensional vector space.  An elliptic boundary condition must cut
it down to a finite dimensional subspace.  Fix $a\in \RR\setminus
\ZZ^{\ge0}$.  Let $\sH_a$ denote the subspace of smooth functions
$u\:\overline\Omega \to\CC$ such that the Fourier expansion
of~$u\res{\partial \Omega }$ has vanishing Fourier coefficient
of~$e^{im\theta }$ if~$m>a$, where we write $z=e^{i\theta }$ on~$\partial
\Omega $.  Then the restriction of~$\dbar$ to~$\sH_a$ has finite dimensional
kernel and cokernel, and it extends to a Fredholm operator on suitable
Sobolev completions.  In other words, restriction to $\sH_a$ is an elliptic
boundary condition.  But because of the Fourier transform in its definition,
it is not local; compare~\eqref{eq:54}.
 
This example generalizes to a Dirac operator~$D_X$ on a compact Riemannian
manifold~$X$ with boundary.  Assume the metric is a product near~$Y=\bX$, and
so decompose
  \begin{equation}\label{eq:90}
     D_X = \sigma \frac{\partial }{\partial t} + D_Y 
  \end{equation}
near the boundary.  Here $t$~is the length coordinate on geodesics normal
to~$Y=\bX$, the algebraic operator~$\sigma $ is Clifford multiplication
by~$dt$, and $D_Y$~is a Dirac operator on~$Y$.  Then the operator $A_Y=\sigma
\inv D_Y$ is self-adjoint.  Let
  \begin{equation}\label{eq:91}
     \bigoplus\limits_{\lambda \in \spec (A_Y)} \!\!\sE_\lambda 
  \end{equation}
be the spectral decomposition of spinors on~$Y$.  For each $a\in \RR\setminus
\spec(A_Y)$, the Atiyah-Patodi-Singer global boundary condition restricts to
the subspace of spinors on~$X$ whose restriction to~$Y$ lies in the
completion of $\oplus _{\lambda <a}\sE_\lambda $.  In fact,
Atiyah-Patodi-Singer choose~$a=0$, as do we in what follows.  If $0\in
\spec(A_Y)$, then one must take into account $\ker D_Y$ separately, as
in~\eqref{eq:92} below.

  \subsection{The Atiyah-Patodi-Singer $\eta $-invariant}\label{subsec:7.3}

With elliptic boundary conditions in hand, Atiyah-Patodi-Singer proceed to
compute $\ind D_X$ for a general Dirac operator on a compact manifold with
boundary.  The problem splits into two pieces: a cylinder near~$\bX$ and
$X\setminus \bX$.  On the cylinder they use the spectral
decomposition~\eqref{eq:91} and the product metric
decomposition~\eqref{eq:90} to convert~$D_X$ to a family of ordinary
differential operators parametrized by $\spec(A_Y)$.  On the complement
of~$\bX$ they use heat kernel methods, as in the local index
theorem~(\S\ref{subsec:6.2}).  Gluing the two regions via a partition of
unity, they prove the following.

  \begin{theorem}[Atiyah-Patodi-Singer~\cite{APS1}]\label{thm:16}
 Let $X$~be a compact Riemannian manifold with boundary, and assume the
Riemannian metric is a product in a neighborhood of~$\bX$.  Let $D_X$ be a
generalized Dirac operator.  Then with respect to the global boundary
conditions, 
  \begin{equation}\label{eq:92}
     \ind D_X = \int_{X}\omega \;-\; \frac{\et{\bX} + h\mstrut _{\bX}}{2}, 
  \end{equation}
where $\omega $~is the Chern-Weil form of the $\Ahat$-genus, $\et{\bX}$~is
the $\eta $-invariant, and $h\mstrut _{\bX}=\dim\ker D_{\bX}$.
  \end{theorem}

\noindent
 To define the $\eta $-invariant on~$Y=\bX$, let $A_Y=\sigma\inv  D_Y$ and form 
  \begin{equation}\label{eq:93}
     \eta (s) = \sum\limits_{\lambda \in \spec(A_Y)\setminus \{0\}}
     (\sign\lambda )\,|\lambda |^{-s},\qquad \Re(s)>\!>0. 
  \end{equation}
This is a Riemannian version of an $L$-function, an echo of the
number-theoretic $L$-function in the signature defect on a Hilbert modular
surface, and a variation on the Riemannian version~\eqref{eq:57} of a $\zeta
$-function.  The infinite sum in~\eqref{eq:93} converges for $\Re(s)>\!>0$,
there is a meromorphic continuation\footnote{Notably, $\eta (s)$~is
holomorphic for $\Re(s)>-1/2$.}  to the complex $s$-line, and $s=0$~is a
regular point.  Define $\et Y=\eta (0)$.  We remark that the main theorem
in~\cite{APS1} applies to more general first-order elliptic differential
operators.
 
The index theorem~\ref{thm:16} simplifies for the signature
operator~\eqref{eq:78}.  Let $X$~be a compact oriented Riemannian
$4k$-manifold with product metric near~$\bX$.  The symmetric bilinear
form~\eqref{eq:8} is nondegenerate restricted to the image of
$H^{2k}(X,\bX;\RR)$ in $H^{2k}(X;\RR)$, and $\Sign(X)$~is its signature.

  \begin{corollary}[Atiyah-Patodi-Singer~\cite{APS1}]\label{thm:17}
 In this situation, 
  \begin{equation}\label{eq:94}
     \Sign(X) = \int_{X}\omega \;-\;\et{\bX}, 
  \end{equation}
where $\omega $~is the Chern-Weil form of the $L$-genus, and $\et{\bX} $~is
the $\eta $-invariant of the self-adjoint operator on $\Omega
^{\textnormal{even}}(\bX)$ given by $(-1)^{k+q+1}(*d-d*)$ on~$\Omega ^q(\bX)$. 
  \end{corollary}

\noindent
 In particular, $-\eta $~is the signature defect~\eqref{eq:87}.  Note that
the signature defect is a \emph{spectral invariant} of a natural differential
operator on~$\bX$, a property which is not apparent from its definition. 

  \begin{remark}[]\label{thm:18}
 The signature defect also plays a star role in two of Atiyah's later
papers~\cite{A8,A9}.  
  \end{remark}

The $\eta $-invariant is our first example of a \emph{geometric} invariant of
a Dirac operator.  To illustrate, recall Example~\ref{thm:24}.  Let $Y\to S$
be a proper fiber bundle of odd relative dimension equipped with a relative
Riemannian spin structure.  From this geometric data we obtain a family of
self-adjoint Dirac operators parametrized by~$S$.  The lowest piece of the
topological index is a \emph{homotopy class} of maps $S\to\RZ$.  The
expression
  \begin{equation}\label{eq:95}
     \xia{Y/S} = \frac{\et{Y/S} + h\mstrut
     _{Y/S}}{2}\!\!\!\!\pmod1\:S\longrightarrow \RZ 
  \end{equation}
from~\eqref{eq:92} refines the homotopy class to a specific map, the
geometric invariant in question.  Also, Theorem~\ref{thm:16} implies that the
differential of~\eqref{eq:95} is
  \begin{equation}\label{eq:96}
     d\xia{Y/S} = \int_{Y/S}\omega , 
  \end{equation}
where $\omega $ is the differential form in~\eqref{eq:92}.  This is a kind of
``curvature'' formula for the geometric invariant~$\xia{Y/S}$; there is an
analogue for other geometric invariants of Dirac operators.

Atiyah-Patodi-Singer~\cite[\S7]{APS3}, in collaboration with Lusztig, gave
another analytic computation of the homotopy class of the map $S\to\RZ$.  For
any loop $\cir\to S$ they prove the winding number of the composite $\cir\to
S\to \RZ$ is the \emph{spectral flow} of the pullback family of Dirac
operators parametrized by~$\cir$.  The spectral flow counts with sign the
integer jumps in the $\xi $-invariant~\eqref{eq:95} as we travel
around~$\cir$.  Alternatively, the union of the spectra of the Dirac
operators is a closed subset $C \subset \cir\times \RR$; the spectral flow is
the intersection number of~$C$ with~$\cir\times \{0\}$.

  \begin{example}[]\label{thm:25}
 The simplest nontrivial spectral flow occurs for the family of complex
self-adjoint Dirac operators 
  \begin{equation}\label{eq:100}
     D_s = \sqrt{-1} \frac{d}{dx} + s,\qquad 0\le s\le 1, 
  \end{equation}
acting on~$\RR/2\pi \ZZ$ with coordinate~$x$.  The operator~$D_1$ is
isomorphic to~$D_0$: conjugate by the multiplication
operator~$e^{\sqrt{-1}x}$.  The union of spectra~ $C\subset \cir\times \RR$
~is a helix.
  \end{example}

  \subsection{Secondary geometric invariants}\label{subsec:7.4}

The most elementary secondary invariant is the total geodesic curvature of a
curve in a Riemannian 2-manifold; it appears in the Gauss-Bonnet
formula~\eqref{eq:85}.  The associated primary topological invariant is the
Euler number.  The generalizations below are more akin to the mod~$\ZZ$
reduction of the total geodesic curvature, which is---up to a sign---the
holonomy of the Levi-Civita connection. 
 
Let $G$~be a Lie group with finitely many components, $\pi \:P\to X$ a
principal $G$-bundle, and $\Theta \in \Omega ^1(P;\mathfrak{g})$ a connection
form.  Let $p\in \bigl(\Sym^k\mathfrak{g}^* \bigr)^G$ be an Ad-invariant
polynomial on the Lie algebra~$\mathfrak{g}$.  As we have already used,
Chern-Weil associate to this data a closed differential form $p(\Omega )\in
\Omega ^{2k}(X)$ that depends only on the curvature~$\Omega $ of~$\Theta $.
Furthermore, it is natural in the connection~$\Theta $.  Its de Rham
cohomology class $\bigl[p(\Omega )\bigr]\in H^{2k}(X;\RR)$ is independent
of~$\Theta $, so is an invariant of the principal bundle~$\pi $.  If $G$~is
compact, which we now assume, then this invariant is derived from a
characteristic class $c_p\in H^{2k}(BG;\RR)$ in the cohomology of the
classifying space of~$G$.  In~1972, Chern-Simons~\cite{CS} introduced a
secondary geometric invariant attached to a refinement of~$c_p$ to an
integral cohomology class $c\in H^{2k}(BG;\ZZ)$.  (Refinements exist only if
the periods of~$c_p$ are integers.)  In this situation, the primary
$\ZZ$-valued invariant is a characteristic number of a principal bundle $\pi
\:P\to X$ over a $2k$-dimensional closed oriented manifold.  The secondary
$\RZ$-valued \emph{Chern-Simons invariant} is defined for $\pi \:Q\to Y$ with
connection~$\Theta $, where $Y$ is a closed oriented $(2k-1)$-dimensional
manifold.  The secondary $\RZ$-valued invariant depends on the connection,
whereas the primary $\ZZ$-valued invariant is topological.

  \begin{example}[]\label{thm:19}
 For~$N\in \ZZ^{\ge3}$ let $p_1\in H^4(\BSO_N;\ZZ)$ be the universal first
Pontrjagin class of a principal $\SO_N$-bundle.  Working intrinsically---that
is, with the tangential geometry of manifolds---the $\ZZ$-valued primary
invariant of a closed oriented 4-manifold~$X$ is
  \begin{equation}\label{eq:97}
     p_1(W)[W]. 
  \end{equation}
The secondary invariant~$\Gamma (Y)$ on a closed oriented Riemannian
3-manifold~$Y$ is the Chern-Simons invariant of its Levi-Civita connection.  In
this case~\cite[\S6]{CS}, $\Gamma (Y)$~is a conformal invariant and an
obstruction to the existence of a conformal immersion $Y\to\EE^4$, where
$\EE^4$~ is Euclidean 4-space.
  \end{example}

  \begin{remark}[]\label{thm:20}
 The Chern-Simons invariant finds a natural expression in \emph{differential
cohomology} \cite{ChS,F2,HS}, which unifies the primary and secondary
invariants in a single framework. 
  \end{remark}

The Atiyah-Patodi-Singer $\eta $-invariant is a secondary invariant analogous
to the Chern-Simons invariant, but in index theory rather than the theory of
characteristic classes.  Let $Y$~be a spin Riemannian manifold and $\pi
\:Q\to Y$ a principal $G$-bundle with connection.  Whereas the Chern-Simons
story begins with an integral cohomology class $c\in H^{\bullet
}(BG;\ZZ)$, to define the $\eta $-invariant we begin with a
complex\footnote{There is a refinement to real representations and real
$KO$-theory.} linear representation $\rho \:G\to\Aut(V)$; the isomorphism
class of~$\rho $ is an element of the equivariant $K$-theory group
$K^0_G(\pt)$.  From the beginning we see the Chern-Simons invariant pertains to
integer cohomology, while the $\eta $-invariant pertains to $K$-theory.  Form
the Dirac operator on~$Y$ coupled to the associated vector bundle $V_Q\to Y$
with its inherited covariant derivative.  Then the $\RZ$-valued invariant 
  \begin{equation}\label{eq:98}
     \xia Y(V) = \frac{\et Y(V) + h\mstrut _Y(V)}{2}\pmod1 
  \end{equation}
depends smoothly on the Riemannian metric and the connection on~$\pi
$.\footnote{For the special case of the signature operator, in which the
kernel---the harmonic forms---have cohomological significance, the
$\RR$-valued invariant is smooth.  (It appears in~\eqref{eq:94}.)}  The
corresponding primary invariant is the $\ZZ$-valued index of the twisted
Dirac operator on manifolds of dimension $\dim Y+1$ equipped with a principal
$G$-bundle.   

  \begin{remark}[]\label{thm:21}
 To obtain a non-topological invariant, we must have $\dim Y$ odd.  The $\eta
$-invariant also leads to new and interesting topological invariants in even
dimensions, for example on unoriented manifolds with a pin
structure~\cite{Gi2}.  
  \end{remark}

  \begin{example}[]\label{thm:22}
 Consider $\dim Y=3$ as in Example~\ref{thm:19}, but now assume $Y$~is a
closed spin Riemannian manifold.  The expression~\eqref{eq:98} for the
standard Dirac operator (no principal $G$-bundle) is the $\RZ$-valued
secondary invariant of the index of the Dirac operator on a closed spin
4-manifold~$X$, which by the index theorem~\eqref{eq:42} is 
  \begin{equation}\label{eq:99}
     \ind D_X = -\frac{1}{24}\,p_1(TW)[W]. 
  \end{equation}
The difference with~\eqref{eq:97} is the rational factor.  The integrality of
the $\Ahat$-genus (recall~\eqref{eq:24}) implies \eqref{eq:99}~is an integer.
Turning to the secondary invariants~$\Gamma (Y)$ and~$\xia Y$, since they are
$\RZ$-valued we cannot multiply~$\Gamma (Y)$ by a non-integral rational
number like~$-1/24$; instead we clear denominators and compare $-24\xia Y$
with~$\Gamma (Y)$.  In the case at hand they agree; a similar comparison in
more general circumstances leads to a spin bordism invariant.  In any case,
we see that the $\eta $-invariant is a more subtle invariant than the
Chern-Simons invariant.  Put differently, the secondary invariants based on
$K$-theory contain refined information over those based on integer
cohomology.  This echoes the stronger topological information derived from
integrality of the primary invariants; see~\S\ref{subsec:2.4}.
  \end{example}

  \begin{remark}[]\label{thm:23}
 \

 \begin{enumerate}

 \item As a concrete illustration of this extra power,
Atiyah-Patodi-Singer~\cite[\S4]{APS2} use the $\eta $-invariant of the
signature operator to refine the Chern-Simons obstruction to conformal
embeddings $Y^3\to\EE^4$.  Also, they show how to use $\eta $-invariants to
construct the Adams $e$-invariant, an invariant of framed bordism.

 \item The primary indices and secondary $\eta $-invariants are unified in
the framework of differential $K$-theory; compare Remark~\ref{thm:20}.
See~\cite{FL} and the references therein.

 \item The Atiyah-Patodi-Singer $\eta $-invariant appears in many contexts in
geometry and beyond.  It also, together with other characters in topological
and geometric index theory, makes many appearances in theoretical physics:
quantum field theory, string theory, and condensed matter theory.  

 \end{enumerate}
  \end{remark}

  \subsection{Determinants of Dirac operators}\label{subsec:7.5}

Before proceeding to Dirac operators, consider a second-order Laplace
operator~$\Delta $ on a closed manifold~$X$, as in~\S\ref{subsec:6.1}.  Then
$\Delta $~has a discrete spectrum consisting of eigenvalues $0\le\lambda
_1\le\lambda _2\le\cdots$, repeated with multiplicity.  Formally, the
determinant of~$\Delta $ is 
  \begin{equation}\label{eq:101}
     \det\Delta \;\textnormal{``$=$''} \;\prod\limits_{m=1}^{\infty}\lambda _m. 
  \end{equation}
(Assume $\lambda _1>0$ or omit the zero eigenvalues to avoid $\det\Delta
=0$.)  Of course, this infinite product diverges.  For example, if $X=\cir$
and $\Delta $~is the usual scalar Laplace operator, then up to a constant the
infinite product is $\prod_{m=1}^{\infty}m^2$ after omitting the zero
eigenvalue.  One way to impart a value to this infinite product, pioneered by
Ray-Singer~\cite{RS} in~1971 and following a technique familiar in complex
analysis~\cite{JL}, is to use the analytic continuation of the $\zeta
$-function~\eqref{eq:57}, which is defined as
  \begin{equation}\label{eq:102}
     \zD(s) = \sum\limits_{m=1}^{\infty}\lambda ^{-s},\qquad s\in \CC, 
  \end{equation}
for $\Re (s)>\!>0$.  Then in the region of absolute convergence of the
infinite sum, we differentiate
  \begin{equation}\label{eq:103}
     -\zeta_{\Delta }'(s) = \sum_{m=1}^{\infty}\lambda ^{-s}\log\lambda , 
  \end{equation}
and then use the regularity of the analytic continuation of~$\zD$ at~$s=0$ to
define
  \begin{equation}\label{eq:104}
     \det\Delta := e^{-\zeta_{\Delta }'(0)} .
  \end{equation}
For the scalar Laplace operator on~$\cir$, the elliptic $\zeta $-function
defined in~\eqref{eq:102} reduces to the Riemann $\zeta $-function, up to a
constant.
 
The first-order Dirac operator~\eqref{eq:41} is not self-adjoint; its domain
and codomain are different.  A finite dimensional model is a linear operator 
  \begin{equation}\label{eq:105}
     T\:V^0\longrightarrow V^1 
  \end{equation}
between different vector spaces~$V^0,V^1$.  There is an induced map
${\textstyle\bigwedge} ^qT\:{\textstyle\bigwedge} ^qV^0\to
{\textstyle\bigwedge} ^qV^1$ on each exterior power.  If $\dim V^0=\dim V^1$,
then the induced map for $q=\dim V^i$ is the determinant 
  \begin{equation}\label{eq:106}
     \det T\:\Det V^0\longrightarrow \Det V^1, 
  \end{equation}
where $\Det V^i = {\textstyle\bigwedge} ^{\dim V^i}(V^i)$ is the determinant
line.  If $V^0=V^1$, then the operator~$\det T$ is multiplication by the
numerical determinant.  But in general\footnote{If $\dim V^0\neq \dim V^1$,
define $\det T=0$.} $\det T$~is an element of a line, namely the
1-dimensional vector space
  \begin{equation}\label{eq:123}
     \Hom(\Det V^0,\Det V^1)
  \end{equation}
called the \emph{determinant line}.  The determinant construction generalizes
to Fredholm operators~\eqref{eq:105}, where now $V^0,V^1$ are typically
infinite dimensional.  The formula~\eqref{eq:123} does not make sense if
$V^0,V^1$ are infinite dimensional; rather, the determinant line is defined
using the finite dimensionality of the kernel and cokernel.  Quillen~\cite{Q}
constructs a determinant line bundle $\pi \:\Det\to \Fred$ over the space of
Fredholm operators together with a continuous section $\det\:\Fred\to \Det$.
It is the next topological invariant of Fredholm operators after the
numerical index~\eqref{eq:34}.

Let $f\:X\to S$ be a proper fiber bundle with fibers of even dimension~$n$,
and assume $f$~is endowed with a relative Riemannian spin structure, as in
Example~\ref{thm:24}.  For a complex vector bundle $V\to X$ with covariant
derivative, the index of the resulting family~$D_{X/S}(V)$ of complex
Clifford linear Dirac operators is computed by the pushforward 
  \begin{equation}\label{eq:107}
     f_!\:K^0(X)\longrightarrow K^{-n}(S). 
  \end{equation}
The numerical index is the image of~$f_![V]$ under $K^{-n}(S)\to H^0(S;\ZZ)$.
Parallel to~\eqref{eq:83}, the next lowest piece of the topological index is
computed by a natural map
  \begin{equation}\label{eq:108}
     \Det\:K^{-n}(S)\longrightarrow H^2(S;\ZZ). 
  \end{equation}
The index theorem for families implies that the \emph{topological}
equivalence class of the Fredholm determinant line bundle
  \begin{equation}\label{eq:109}
     \Det D_{X/S}(V)\longrightarrow S 
  \end{equation}
is $\Det f_![V]\in H^2(S;\ZZ)$. 
 
There is a \emph{geometric} refinement of this topological piece of the
index.  Its analytic expression is a hermitian connection and compatible
covariant derivative on the determinant line bundle~\eqref{eq:109}.  The
\emph{Quillen metric}~\cite{Q} is constructed using the $\zeta $-function
determinant~\eqref{eq:104}, and the covariant derivative uses a similar---but
somewhat more subtle---$\zeta $-function~\cite{BiF1}.  The isomorphism class
of a line bundle with covariant derivative is determined by its holonomy
around loops.  Base change along a loop $\varphi \:\cir\to S$ gives rise to a
fiber bundle $X_\varphi \to \cir$ and vector bundle $V_\varphi \to X_{\varphi
}$.  The bounding spin structure on~$\cir$ combines with the relative spin
structure on~$X_\varphi /\cir$ to produce a spin structure on~$X_\varphi $.
Choose an arbitrary metric~$g\mstrut _{\cir}$ on~$\cir$, and let $X_\varphi
(\epsilon )$ be the manifold~$X_\varphi $ with Riemannian metric $g\mstrut
_{\cir}\,/\,\epsilon ^2\,\oplus \,g\mstrut _{X_\varphi /\cir}$.  (The direct
sum is with respect to the horizontal distribution in the relative Riemannian
structure.)  Then the holonomy around the loop~$\varphi $ is~\cite{BiF2}
  \begin{equation}\label{eq:110}
     \hol_\varphi \Det D_{X/S}(V) = \lim\limits_{\epsilon \to0}e^{-2\pi i\xia
     {X_\varphi (\epsilon )}(V)}, 
  \end{equation}
where the exponent is the Atiyah-Patodi-Singer $\eta
$-invariant~\eqref{eq:98}.  This holonomy formula was inspired by Witten's
global anomaly~\cite{Wi3}.  The curvature of the determinant line bundle is 
  \begin{equation}\label{eq:111}
     \curv\Det D_{X/S}(V) = 2\pi i\int_{X/S}\omega , 
  \end{equation}
where $\omega $~is the Chern-Weil form that represents
$\Ahat(X/S)\smallsmile\ch(V)$; compare~\eqref{eq:96}. 

  \begin{remark}[]\label{thm:26}
 \ 
 \begin{enumerate}

 \item The \emph{adiabatic limit} in~\eqref{eq:110} was introduced
in~\cite{Wi3}.  Other geometric interpretations of Witten's global anomaly
formula were given in~\cite{Che,S2}.

 \item The holonomy formula~\eqref{eq:110} may be regarded as a Fubini
theorem relating the geometric invariants~$\eta $ and~$\Det$, once one knows
that the holonomy of a line bundle $L\to \cir$ is the exponentiated $\eta
$-invariant of $D_{\cir}(L)$.

 \item As in Remark~\ref{thm:23}(2), the isomorphism class of the determinant
line bundle with its metric and covariant derivative can be computed by a
pushforward in differential $K$-theory.  This refines the topological index
theorem which computes the topological isomorphism class of~\eqref{eq:109} as
$\Det f_![V]$.

 \end{enumerate}
  \end{remark}

   \section{Anomalies and index theory}\label{sec:8}

In part inspired by Is Singer's advocacy of theoretical physics as a fertile
ground for geometers, beginning in the late~1970's Michael Atiyah turned his
attention to geometric problems in quantum field theory and, later, string
theory.  Simon Donaldson's paper~\cite{Do} in this volume covers the burst of
activity in the late~1970's and early~ 1980's emanating from the Yang-Mills
equations.  The Atiyah-Singer index theorem is a part of that story, but I
will restrict my exposition here to Atiyah's work in the late 1980's on
anomalies and on his axiomatization of topological field theory.
 
We begin in~\S\ref{subsec:8.1} by arguing that anomalies are an expression of
the projective (as opposed to linear) nature of quantum theory.  The
geometrical and topological link between anomalies and index theory was
initiated in a joint paper of Atiyah-Singer, which we summarize
in~\S\ref{subsec:8.3}.  Atiyah's axioms for topological field theory, and
their relationship to bordism in algebraic topology, are the subjects
of~\S\ref{subsec:8.2}.  We conclude in~\S\ref{subsec:8.4} by mentioning the
modern point of view on anomalies, which ties together the two aforementioned
Atiyah works.  The specialization to spinor fields brings in topological and
geometric index theorems as well, all synthesized in a general anomaly
formula which draws on many of the Atiyah papers we have discussed.

  \subsection{Projectivity and symmetries in quantum mechanics}\label{subsec:8.1}

Anomalies are often said to be the failure of a classical symmetry to hold in
a corresponding quantum system.  More precisely, a quantum symmetry is
projective and the anomaly is the obstruction to linearization.  Quantum
theory is inherently projective, and the anomaly---in a more general sense
than a notion tied to symmetry---encodes the projectivity of a quantum
system.

A quantum mechanical system is specified\footnote{There is a more general
framework for quantum theory using $C^*$-algebras, but for this exposition
the simpler context suffices.} by a triple of data~$(\sP,p,H)$.  The
space~$\sP$ is a projective space, the projectivization~$\PH$ of a complex
separable Hilbert space~$\sH$, but $\sH$~is not singled out.  One way to
define~$\sP$ is to fix~$\sH_0$ and let $\cP$~be the groupoid whose objects
are pairs~$(\sH,\theta )$ of a Hilbert space and a projective linear
isomorphism $\theta \:\PH_0\to \PH$.  Morphisms $(\sH,\theta )\to
(\sH',\theta ')$ are linear isometries $\sH\to \sH'$ whose projectivization
commutes with~$\theta ,\theta '$; they form a torsor over the unitary group
$\TT\subset \CC^{\times }$ of unit norm scalars.  The basepoint~$\sH_0$ is
\emph{not} part of the structure.  Define $\sP$ as the limit of~$\PH$ over
all $(\sH,\theta )\in \cP$.  Also, define the $*$-algebra~$\AP$ (of
observables) as the limit of~$\End\sH$ over $(\sH,\theta )\in \cP$.
($\End\sH$~is the algebra of bounded linear operators on~$\sH$; for
simplicity, we omit unbounded operators from this exposition.)  The
projective space~$\sP$ is the space of \emph{pure states} of the quantum
system.  Embed $\sP\hookrightarrow \AP$ as rank one orthogonal projections.
Then mixed states are convex combinations of pure states.  The function
  \begin{equation}\label{eq:112}
     \begin{aligned} p\:\sP\times \sP&\longrightarrow \;\;[0,1] \\
      L_1,L_2&\longmapsto |\langle \psi _1,\psi _2
     \rangle|^2,\qquad \psi _i\in L_i\end{aligned}  
  \end{equation}
is used to compute transition probabilities.  To define~$p$, choose
$(\sH,\theta )\in \cP$ and identify~$\sP$ with~$\PH$; then $L_i\subset \sH$
is a line and $\psi _i\in L_i$ is a unit norm vector.  The
\emph{Hamiltonian}~$H$ is a self-adjoint element of~$\AP$, usually assumed to
have spectrum bounded below.
 
Fix $(\sH,\theta )\in \cP$.  Let $G(\sH)$~be the infinite-dimensional Lie
group with identity component the group~$U(\sH)$ of unitary automorphisms
of~$\sH$ and off-identity component the torsor of antiunitary automorphisms.
A basic theorem of Wigner asserts that the sequence of Lie group
homomorphisms
  \begin{equation}\label{eq:113}
     1\longrightarrow \TT\longrightarrow
     G(\sH)\xrightarrow{\;\;q\;\;}\Aut(\PH,p)\longrightarrow 1 
  \end{equation}
is a group extension: $q$~is surjective.  Let $G$~be a Lie group of
symmetries of~$(\sP,p)$, i.e., a homomorphism $G\to \Aut(\sP,p)$.  By
pullback, we obtain a group extension
  \begin{equation}\label{eq:114}
     1\longrightarrow \TT\longrightarrow G_{\sH}\longrightarrow G\longrightarrow
     1 
  \end{equation}
and a $\zt$-grading $\epsilon \:G\to \zt$.  The $\zt$-grading is independent
of~$(\sH,\theta )$, as is the isomorphism class of the
extension~\eqref{eq:114}.  The $\zt$-graded group extension~\eqref{eq:114}
can be called the \emph{anomaly}; it measures the projectivity of the
symmetry.

  \begin{remark}[]\label{thm:27}
 \

 \begin{enumerate}

 \item If $(\sP,p,H)$ has a classical limit, and the symmetry persists in the
limit, then it is the Lie group~$G$ which acts on the classical system.  The
passage from~$G$ to~$G_\sH$ is what was referenced in the first paragraph of
this section.

 \item Suppose given a family of quantum systems with parameter manifold~$S$,
so in particular a fiber bundle $\sP\to S$ of projective Hilbert spaces.  In
some physical situations one wants to ``integrate over~$S$'' to form a new
quantum mechanical system.  If we write $\sP\to S$ as the projectivization of
a vector bundle $\sE\to S$ of Hilbert spaces, then the state space of the
integrated system is the space of $L^2$~sections of $\sE\to S$, assuming a
measure on~$S$.  The anomaly is the obstruction to lifting $\sP\to S$ to a
vector bundle, a necessary first step to integrate out~$S$.  On the other
hand, in the absence of integrating out~$S$ the anomaly is not an obstruction
but rather a feature of a family of quantum systems, useful in many
contexts. 

 \item The isomorphism class of the central extension~\eqref{eq:114} lives in
the cohomology group $H^2(G;\TT)$.  (The type of cohomology depends on the
type of group.  Nothing is lost here by assuming that $G$~is a finite group.)
If we drop a cohomological degree, then $H^1(G;\TT)$ is the group of
1-dimensional unitary representations.  Therefore, the projectivity measured
in~\eqref{eq:114} is obtained via a 2-step procedure starting with linear
actions of~$G$ on~$\sH$: (1)~replace~$\sH$ by a 1-dimensional vector space,
and (2)~interpret 1-dimensional representations of~$G$ cohomologically and
raise the cohomological degree by~1.  We will see an analogous procedure in
quantum field ~(\S\ref{subsec:8.4}).

 \end{enumerate}
  \end{remark}

Let $t_0<t_1<\cdots< t_{n+1}$ be real numbers, thought of as points on the
affine time line, and let $A_1,\dots ,A_n\in \AP$.  Fix initial and final
pure states $L_0,L_{n+1}\in \sP$.  A basic quantity of interest in quantum
mechanics is the probability\footnote{If the linear operator $\evo{n+1}n
A_n\cdots\evo21 A_1\evo10$ is zero on~$L_0$, the probability is zero.}
  \begin{equation}\label{eq:115}
     p\bigl(L_{n+1}\,,\, \evo{n+1}n A_n\cdots\evo21 A_1\evo10 L_0 \bigr),
  \end{equation}
where $\hbar$~is Planck's constant.  Physical questions may be phrased in
terms of these probabilities.  Fix $(\sH,\theta )\in \cP$, identify
$\sP\approx \PH$, and choose unit norm vectors $\psi _0\in L_0$ and
$\psi _{n+1}\in L_{n+1}$.  The \emph{amplitude}, or \emph{correlation
function}, of this data is the complex number
  \begin{equation}\label{eq:116}
     \bigl\langle \psi _{n+1}\,,\, \evo{n+1}n A_n\cdots\evo21 A_1\evo10 \psi _0
     \bigr\rangle_{\sH} .
  \end{equation}
More invariantly, the data $L_0,L_{n+1},t_0,\dots ,t_{n+1},A_1,\dots ,A_n$
determine a hermitian line~$\sL$, and the amplitude is an element of~$\sL$
whose norm is the probability~\eqref{eq:115}.

  \begin{remark}[]\label{thm:28}
 The fact that the amplitudes of~$(\sP,p,H)$ lie in a line~$\sL$ without a
distinguished basis element is another aspect of the projectivity of quantum
mechanics.  The lines~$\sL$ form a line bundle over a parameter space of
data.  These line bundles are part of the anomaly of the quantum mechanical
system, a counterpart for correlation functions of the projective bundles in
Remark~\ref{thm:27}(2).
  \end{remark}

 The context for anomalies sketched in this section is soft, much as is
index theory for general Fredholm operators.  Quantum field theory brings in
the geometry of Wick-rotated spacetimes, and in that context anomalies also
exhibit more geometry, much as does index theory for Dirac operators.  This
is more than an analogy in the case of fermionic fields, to which we now
turn.

  \subsection{Spinor fields and anomalies in quantum field theory}\label{subsec:8.3}

Geometric links between anomalies for spinor fields and index theory were
forged in a 1984 paper of Atiyah-Singer~\cite{AS7}.  (Other contemporaneous
papers, such as \cite{AgW,AgG,Lt2} also brought index polynomials into the
theory of anomalies.)  The setup is gauge theory in physics.  (By that time
Atiyah had already achieved many results in mathematical gauge theory; see
the article~\cite{Do} by Donaldson in this volume.)  We first summarize
their work and then relate it to the geometric picture of anomalies.
 
Let $n$~be an even positive integer and $X$~a closed Riemannian spin manifold
of dimension~$n$.  Suppose $P\to X$ is a principal bundle with structure
group a compact Lie group~$G$.  Let $\sA$~be the infinite dimensional affine
space of connections on $P\to X$, and let $\sG$~be the group of gauge
transformations: automorphisms of $P\to X$ which act as the identity on~$X$.
Then $\sG$~acts on~$\sA$.  Assume $X$~is connected, fix a
basepoint~$p_0\in P$, and let $\sG_0\subset \sG$ be the subgroup of gauge
transformations which fix~$p_0$.  Then $\sG_0$~acts freely on~$\sA$, and in
the sequence 
  \begin{equation}\label{eq:117}
     \frac{\sA\times P}{\sG_0}\xrightarrow{\;\;\pi
     \;\;}\frac{\sA}{\sG_0}\times X
     \xrightarrow{\;\;\pr_1\;\;}\frac{\sA}{\sG_0} ,
  \end{equation}
$\pi $~is a principal $G$-bundle and $\pr_1$~is a product fiber bundle with
fiber~$X$.  To a unitary representation of~ $G$ we associate a vector bundle
to~$\pi $.  This data determines a topological index
  \begin{equation}\label{eq:118}
     \ind\in K^{-2n}(\sA/\sG_0)\cong K^0(\sA/\sG_0). 
  \end{equation}

Atiyah-Singer construct closed differential forms which represent
$\ch(\ind)\in H^{\bullet }(\sA/\sG_0;\RR)$ as follows.  Fix a bi-invariant
Riemannian metric on~$G$.  Then for each connection~$A\in \sA$, the
manifold~$P$ has a Riemannian metric which makes the projection $P\to X$ a
Riemannian submersion in which the horizontal subspaces of the connection are
orthogonal to the fibers.  Use the $L^2$~metric on~$\sA$ to form a warped
product metric on~$\sA\times P$.  Then $\sG_0\times G$ acts by isometries,
and a connection on~$\pi $ results by taking orthogonals to the $G$-orbits on
the quotient by~$\sG_0$.  The associated vector bundle inherits a covariant
derivative, so there is a family of Dirac operators on~$X$ parametrized
by~$\sA/\sG_0$, and the index theorem for families (\S\ref{subsec:5.4})
implies that ind in~\eqref{eq:118} equals its analytic index.  The Chern-Weil
procedure produces the desired differential forms from the differential
geometric data.  In particular, $\ch_1(\ind)\in H^2(\sA/\sG_0;\RR)$ is
represented by
  \begin{equation}\label{eq:119}
     \int_{X}\omega 
  \end{equation}
for a $(2n+2)$-form~$\omega $ on $(\sA/\sG_0)\times X$.  Atiyah-Singer
transgress~$\omega $ to a $(2n+1)$-form on $\sG_0\times X$, and it is this
differential form which appears in the contemporaneous physics literature as
the anomaly of a spinor field in quantum field theory.  They also relate
their result to determinants.  In particular, they interpret the cohomology
class $c_1(\ind)\in H^2(\sA/\sG_0;\ZZ)$ of~\eqref{eq:119} as the isomorphism
class of the determinant line bundle of the family of Dirac operators on~$X$
parametrized by~$\sA/\sG_0$.
 
The determinant, or more generally pfaffian, of a Dirac operator arises
directly in Wick-rotated quantum field theories with spinor fields, as we now
sketch in a general context.  Let $n$~be a positive integer, $G$~a compact
Lie group, and $\rho $~a representation of~$G$.  Suppose
  \begin{equation}\label{eq:120}
     P\xrightarrow{\;\;\pi \;\;}X\xrightarrow{\;\;f\;\;}S 
  \end{equation}
 is a principal $G$-bundle~$\pi $ with connection and a proper fiber
bundle~$f$ of relative dimension~$n$, the latter equipped with a relative
Riemannian spin structure.  As in~\S\ref{subsec:5.4}, this data produces a
family of Dirac operators~$D_{X/S}$ parametrized by~$S$.  The Feynman
integral over the spinor field~$\psi $ is an infinite dimensional variant of
a standard Gaussian integral:
  \begin{equation}\label{eq:121}
     \int_{\sS_s} e^{-(\psi ,D_s\psi )}. 
  \end{equation}
For each~$s\in S$ we are meant to integrate over the infinite dimensional
vector space~$\sS_s$ of spinor fields.  The Dirac operator~$D_s$ appears in
the exponent, which is a skew-symmetric bilinear form
on~$\sS_s$.\footnote{Physicists begin with a relativistic theory on Minkowski
spacetime in Lorentz signature, and then they Wick rotate to Riemannian
manifolds.  What emerges naturally is the bilinear form in the exponent.}  By
analogy with finite dimensional integrals, \eqref{eq:121}~is defined to be
$\pfaff D_s$, the pfaffian of the Dirac operator.  In special cases, such as
the Atiyah-Singer setup above, the pfaffian reduces to a determinant, and we
make that simplification here.  Then the result of the
integral~\eqref{eq:121} is the section $\det D_{X/S}$ of the determinant line
bundle  
  \begin{equation}\label{eq:122}
     \Det D_{X/S}\longrightarrow S, 
  \end{equation}
as in~\S\ref{subsec:7.5}. 

  \begin{remark}[]\label{thm:29}
 \

 \begin{enumerate}

 \item The integral~\eqref{eq:121} is an example of a correlation function in
Feynman's approach to quantum field theory.  The fact that it is an element
of a complex line, rather than a complex number, is precisely the situation
of anomalous amplitudes discussed in~\S\ref{subsec:8.1}.  Here the
determinant line bundle~\eqref{eq:122}, including its metric and covariant
derivative, is the anomaly.

 \item The universal parameter space~$S$, for a fixed manifold~$X$, is the
space of Riemannian metrics and $G$-connections on~$X$.  In a quantum field
theory one often wants to integrate over the metric or connections or both.
The anomaly is an obstruction to doing so.

 \end{enumerate}
  \end{remark}

The anomaly is the obstruction to a trivialization of~\eqref{eq:122}.  The
topological first Chern class of the index discussed by Atiyah-Singer
(cf.~\cite{S3}) obstructs a nonzero section, but one needs something
sharper.  In~\cite{F3} the anomaly is identified as the obstruction to a
\emph{flat} section, relative to the natural connection~\cite{BiF1}.  This
leads to an interpretation of the 2-form~\eqref{eq:119} computed by
Atiyah-Singer as the curvature of this natural connection.  Physicists call
this local obstruction to a flat section the \emph{local anomaly}; the
holonomy is the \emph{global anomaly} of Witten~\cite{Wi3}.  Indeed, the
developments in geometric index theory recounted in~\S\ref{subsec:7.5} were
directly inspired by this geometry of anomalies.

  \begin{remark}[]\label{thm:33}
 The anomaly as an obstruction to lifting a bundle of projective spaces
to a vector bundle (Remark~\ref{thm:27}(2)) is measured by the next invariant
in geometric index theory after the determinant line bundle---the Dirac
gerbe~\cite{Lt1,Bu}.  We refer the reader to~\cite{FS,Seg2,NAg}. 
  \end{remark}

  \subsection{Topological field theory and bordism}\label{subsec:8.2}

Atiyah's engagement with quantum field theory went well beyond the index
theorem.  One particularly influential paper~\cite{A11} sets out axioms for
\emph{topological} quantum field theories, parallel to axioms introduced
previously by Segal~\cite{Seg3} for 2-dimensional conformal field theories.
One key impetus was Witten's quantum Chern-Simons theory~\cite{Wi4}, which
places the Jones invariants of knots in a manifestly 3-dimensional framework.
But there were many other examples too that Atiyah abstracted into his
axioms.  His paper~\cite{A11} is dedicated to Thom, and indeed bordism theory
is very much at the forefront.  Atiyah~\cite{A12} writes the following about
this paper:

 \begin{quote} 
 Because mathematicians are frightened by the Feynman integral and are
unfamiliar with all the jargon of physicists there seemed to me to be a need
to explain to mathematicians what a topological quantum field theory really
was, in user-friendly terms.  I gave a simple axiomatic treatment (something
mathematicians love) and listed the examples that arise from physics.  The
task of the mathematician is then to construct, by any method possible, a
theory that fits the axioms.  I like to think of this as analogous to the
Eilenberg-Steenrod axioms of cohomology, where one can use simplicial, Cech
or de Rham methods to construct the theory.  This last is closest to physics
but the others have some advantages.  In the quantum field theory context,
where things are vastly more difficult, the combinatorial approach is so far
the only one that has been made to work (for the Jones polynomials).
 \end{quote}

\noindent
 Here we give a concise version of the axioms.
 
As motivation, recall the signature of a closed oriented manifold of
dimension~$4k$ for some $k\in \ZZ^{\ge0}$.  As used crucially in the proof of
Theorem~\ref{thm:3}, the signature is a \emph{bordism invariant}, that is, a
homomorphism of abelian groups
  \begin{equation}\label{eq:124}
     \Sign\:\Omega _{4k}(\SO)\longrightarrow \ZZ, 
  \end{equation}
where $\Omega _{4k}(\SO)$ is Thom's bordism group of closed oriented
$4k$-manifolds.  A topological field theory is a ``categorified bordism
invariant''.  Fix a nonnegative integer~$d$.  (The relation to~$n$
in~\S\ref{subsec:8.3} is $n=d+1$.)  Let $\Bd$ be the following category,
first introduced by Milnor~\cite{Mi2}.  The objects are closed $d$-manifolds.
If $Y_0,Y_1$ are two such, then a morphism $Y_0\to Y_1$ is represented by a
compact $(d+1)$-manifold~$X$ with boundary partitioned as $\partial X=Y_0\amalg
Y_1$.  In other words, $X$~is a bordism from~$Y_0$ to~$Y_1$.  Diffeomorphic
bordisms rel boundary represent the same morphism.  Composition glues
bordisms and disjoint union of manifolds provides a symmetric monoidal
structure.

  \begin{remark}[]\label{thm:30}
 \

 \begin{enumerate}

 \item If we declare objects $Y_0,Y_1$ of~$\Bd$ to be equivalent if there
exists a morphism $Y_0\to Y_1$, then the set of equivalence classes is the
bordism group~$\Omega _d$.  In this sense, $\Bd$ ``categorifies''~$\Omega
_d$.

 \item A small variation yields bordism categories with tangential structure,
such as an orientation.

 \end{enumerate}
  \end{remark}

We can now state the axioms.  Let $k$~be a field and let $\Vk$~be the
category of $k$-vector spaces and linear maps.  (In quantum theories,
$k=\CC$.)  The operation of tensor product defines a symmetric monoidal
structure on~$\Vk$.

  \begin{definition}[Atiyah~\cite{A11}]\label{thm:31}
 A \emph{topological field theory} is a symmetric monoidal functor 
  \begin{equation}\label{eq:125}
     F\:\Bd\longrightarrow \Vk. 
  \end{equation}
  \end{definition}

\noindent
 This definition is sometimes referred to as the \emph{Atiyah-Segal Axiom
System}, and with suitable modifications and extensions it is believed to
apply widely to field theories in mathematics and physics.  The viewpoint in
these axioms is very different from what one sees in physics texts.  The
Atiyah-Segal Axiom System has provided a generation of mathematicians with a
point of entry to this physics, they are the structure upon which many
mathematical developments have been built, and they have illuminated
geometric aspects of quantum field theories in physics as well.
 
To illustrate Definition~\ref{thm:31}, we show how to extract numerical
invariants of a normally framed knot $K\subset M$ in a closed 3-manifold~$M$
from a field theory~$F$ with~$d=2$.  Let $X$~be the 3-manifold obtained
from~$M$ by removing an open tubular neighborhood of~$K$.  The result is a
bordism $X\:\partial X\to \emptyset ^2$ from~$\partial X$ to the empty
2-manifold.  A normal framing of~$K$ provides an isotopy class of
diffeomorphisms $\partial X\approx\cir\times \cir$.  Hence the value of the
field theory~$F$ on~$X$ is a linear map
  \begin{equation}\label{eq:126}
     F(X)\:V\longrightarrow k, 
  \end{equation}
where $V=F(\cir\times \cir)$ is the vector space attached to the standard
2-torus.  For each vector~$\xi \in V$, which may be viewed as a ``label''
attached to~$K$, we obtain a numerical invariant~$F(X)(\xi )$.  The Jones
invariants of knots are of this type.  This is one of the key observations
in~\cite{Wi4}.

  \subsection{Synthesis}\label{subsec:8.4}

We conclude by bringing together the Atiyah-Singer work on anomalies
(\S\ref{subsec:8.3}), the Atiyah-Segal Axiom System for quantum field
theory (\S\ref{subsec:8.2}), and the index theorems in~\S\S\ref{sec:4},
\ref{sec:5}, \ref{sec:7}.   
 
One starting point is Remark~\ref{thm:29}(2), which tells that an anomaly
must be trivialized to construct a quantum field theory by integrating over
certain fields, such as metrics or connections.  (Such integrals are
problematic analytically, but the anomaly and trivializations are
mathematically well-defined.)  Now if the resulting quantum field theory is
to be \emph{local}---and locality is a characteristic feature of quantum
field theories, then the trivializations of the anomaly must be coherent in
the background data~\eqref{eq:120}.  That coherence is precisely what is
expressed in the Atiyah-Segal Axiom System and its extensions.  This is one
line of reasoning which leads to the realization that an anomaly itself is a
quantum field theory,\footnote{This is not quite universally true: The
anomaly of an $n$-dimensional field theory may not be a full
$(n+1)$-dimensional theory, but may only be defined on manifolds of
dimension~$\le n$.}  albeit of a very special type.  For the spinor field
in~\S\ref{subsec:8.3}, the determinant lines in the fibers of~\eqref{eq:122}
are 1-dimensional state spaces in an $(n+1)$-dimensional field theory.  This
\emph{anomaly theory} is \emph{invertible}, but is not necessarily
topological.  (An invertible field theory~\eqref{eq:125} factors through the
subgroupoid $\Line_k\subset \Vk$ of lines and invertible linear maps.)

  \begin{remark}[]\label{thm:34}
 We arrive at the same picture by following the ideas of~\S\ref{subsec:8.1}.
Namely, a field theory in the form~\eqref{eq:125} is a \emph{linear}
representation of bordism, but quantum theory is \emph{projective} and the
anomaly measures the projectivity.  Furthermore, Remark~\ref{thm:27}(3) gives
a roadmap to locate this measurement.  Here the bordism category plays the
role of the group~$G$.  First, replace~$\Vk$ by the category of 1-dimensional
vector spaces and invertible linear maps.  What results is an
\emph{invertible} field theory.  Second, we interpret an invertible field
theory cohomologically and raise the cohomological degree by~1.  The
cohomological interpretation was introduced in~\cite{FHT};
see~\cite[\S5]{FH}.  We arrive at the same conclusion: the anomaly, or
measurement of projectivity, of a field theory is an invertible
$(n+1)$-dimensional field theory.\footnote{More precisely, an $n$-dimensional
field theory is a representation of a bordism category~$\sB$.  Following the
logic of this paragraph, the anomaly is defined on the \emph{same} bordism
category~$\sB$, so it is a ``once-categorified invertible $n$-dimensional
field theory'', as indicated in the previous footnote.  Typically, the
anomaly extends to a full invertible $(n+1)$-dimensional field theory, but
that is not required.}

  \end{remark}

This already brings together the aforementioned 1980s work of Atiyah, but we
can go much further.  Whereas a general field theory is a functor between
symmetric monoidal categories, an invertible field theory can be formulated
in stable homotopy theory, as a map of spectra.  The domain, rather than a
bordism category, is a bordism spectrum of the type\footnote{This statement
is for unitary theories.} introduced by Thom.  (We remark that
Atiyah~\cite{A10} put bordism and cobordism in the context of generalized
homology theories.)  For spinor fields the domain spectrum is $M\!\Spin$ or a
close variant.  The general formula for the anomaly of a spinor
field~\cite[Conjecture~9.70]{FH}---conjectural as a mathematical assertion
until more foundations are laid---brings in the Atiyah-Bott-Shapiro map
$M\!\Spin\to KO$, as well as all of the aforementioned ingredients.  Implicit
in it are the various topological and geometric index invariants and index
theorems that we have surveyed in this article.


 \bigskip\bigskip
 
\providecommand{\bysame}{\leavevmode\hbox to3em{\hrulefill}\thinspace}
\providecommand{\MR}{\relax\ifhmode\unskip\space\fi MR }
\providecommand{\MRhref}[2]{%
  \href{http://www.ams.org/mathscinet-getitem?mr=#1}{#2}
}
\providecommand{\href}[2]{#2}

  \end{document}